\documentclass[onefignum,onetabnum]{siamonline250106}

\usepackage{graphicx} 
\usepackage{subfigure}
\usepackage{enumitem}   
\usepackage{amsmath,amssymb}
\usepackage{blkarray, bigstrut, bigdelim}
\usepackage{mathtools,mathrsfs}
\usepackage{multirow}
\usepackage{siunitx}
\usepackage{xspace}
\usepackage{comment}
\usepackage{enumitem}
\usepackage{tikz}
\usetikzlibrary{matrix,calc}
\usepackage{algorithm}%
\usepackage[noend]{algpseudocode}%
\algrenewcommand\textproc{}

\newcommand{\eg}{e.\,g.\xspace}
\newcommand{\ie}{i.\,e.\xspace}
\newcommand{\norm}[2]{\left\Vert {#2} \right\Vert_{#1}}

\newcommand{\given}{\left.\hspace*{-0.0833em}\middle|\hspace*{-0.0833em}\right.}

\newcommand{\bvec}[1]{\textbf{#1}}
\newcommand{\bvecS}[1]{\ensuremath{\boldsymbol{#1}}}
\newcommand{\design}{\bvec{e}} 
\newcommand{\params}{\bvec{m}} 
\newcommand{\refparam}{z} 
\newcommand{\refparams}{\bvec{\refparam}} 
\newcommand{\unitvec}{\bvec{u}} %
\newcommand{\targetparam}{x}
\newcommand{\targetparams}{\bvec{\targetparam}} 
\newcommand{\ratioparams}{\bvec{w}} 
\newcommand{\data}{\bvec{y}} 
\newcommand{\noise}{\bvecS{\eta}} 
\newcommand{\m}{\mathrm{m}}
\newcommand{\LISDF}{\bvec{U}} 
\newcommand{\LISDD}{\bvec{V}} 
\newcommand{\source}{S} 

\newcommand{\Ne}{n_{e}} 
\newcommand{\Ns}{s} 
\newcommand{\Nd}{n_{y}} 
\newcommand{\Nm}{n_{m}} 
\newcommand{\TTRank}{R} 
\newcommand{\LISdim}{r} 
\newcommand{\paramSpace}{\mathbb{R}^{\Nm}}
\newcommand{\paramSubspace}{\mathcal{M}}
\newcommand{\designSpace}{\mathcal{E}}
\newcommand{\dataSpace}{\mathcal{Y}}
\newcommand{\targetSpace}{\mathcal{X}}

\newcommand{\PtO}{\mathcal{G}} 
\newcommand{\discPtO}{\bvec{G}} 
\newcommand{\I}[1]{\bvec{I}_{#1}}

\newcommand{\Proj}{\bvec{P}} 
\newcommand{\TTCore}[2]{\bvec{F}_{#1}(#2)}
\newcommand{\pushforward}[1]{\mathop{#1_\sharp}}
\let\pullback\undefined
\newcommand{\pullback}[1]{\mathop{#1^\sharp}}
\newcommand{\cS}{\mathcal{S}} 
\newcommand{\cT}{\mathcal{T}} 
\newcommand{\cQ}{\mathcal{Q}} 
\newcommand{\Jac}{\nabla} 

\newcommand{\refD}{\rho} 
\newcommand{\target}{\pi} 
\newcommand{\targetUnnorm}{\gamma} 
\newcommand{\likelihood}[2]{\mathcal{L}(#1\given#2)}
\newcommand{\Rlikelihood}[2]{\widetilde{\mathcal{L}}(#1\given#2)}
\newcommand{\ratio}[1]{q^{#1}} 
\newcommand{\Cnoise}{\bvecS{\Gamma}_{\noise}} 
\newcommand{\Cnoiseinv}[1]{{\Cnoise(#1)}^{-1}} 

\newcommand{\Cpr}{\mathcal{C}_{0}}

\newcommand{\EIG}{\Psi} 
\newcommand{\distKL}[2]{\mathcal{D}_\textup{KL}{(}{#1}{\|}{#2}{)}}
\newcommand{\distH}[2]{\mathcal{D}_\textup{H}{(}{#1,\,#2}{)}}
\newcommand{\Hist}{\bvec{H}} 
\newcommand{\Hmat}{\mathcal{H}} 
\newcommand{\Expect}[2]{\mathbb{E}_{#1}{[}{#2}{]}}

\DeclareMathOperator*{\argmax}{arg\,max}
\DeclareMathOperator*{\Argmax}{Arg\,max}

\DeclareMathOperator{\Cov}{Cov}

\title{Subspace accelerated measure transport methods for fast and scalable sequential experimental design, with application to photoacoustic imaging\footnote[0]{\funding{The work of RH, KK, and RS was partially funded by the Carl Zeiss Stiftung through the project ``Model-Based AI: Physical Models and Deep Learning for Imaging and Cancer Treatment''. The work of TC is partially supported by the Australian Research Council grant FT250100199 and the Romberg Visiting Scholar Programme of Heidelberg University.}}}

\def\addressC{Interdisciplinary Center for Scientific Computing (IWR), Heidelberg University, Heidelberg, Germany}
\def\addressE{Institute for Mathematics, Heidelberg University, Heidelberg, Germany}
\def\addressD{School of Mathematics and Statistics, University of Sydney, Australia.}

\author{Tiangang Cui\footnotemark[1]
    \and Karina Koval\footnotemark[2]
     \and Roland Herzog\footnotemark[2] \footnotemark[3]
    \and Robert Scheichl\footnotemark[3] \footnotemark[2]}

\footnotetext[2]{\addressC}
\footnotetext[1]{\addressD}
\footnotetext[3]{\addressE}

\begin{document}

\allowdisplaybreaks

\maketitle

\begin{abstract}
We propose a novel approach for sequential optimal experimental design (sOED) for Bayesian inverse problems involving expensive models with high-dimensional unknown parameters. This work focuses on designs that maximize the expected information gain (EIG) from prior to posterior, a task that is computationally very challenging in non-Gaussian settings. This challenge is amplified in sOED, as the incremental expected information gain (iEIG) must be repeatedly approximated across distinct stages, with both prior and posterior distributions being intractable. To address this, we derive a general-purpose, derivative-based upper bound for the iEIG, which not only guides design placement but also enables the construction of projectors onto likelihood-informed subspaces, facilitating parameter dimension reduction. By combining this approach with conditional measure transport maps for the sequence of posteriors, we develop a unified sOED and amortized inference framework scalable to high- and infinite-dimensional problems. Numerical experiments for two inverse problems governed by partial differential equations (PDEs) demonstrate the effectiveness of designs by maximizing the proposed bound.
\end{abstract}



\section{Introduction}

We consider the sequential design of a finite number of experiments for Bayesian inverse problems.
In sequential optimal experimental design (sOED), experimental conditions are chosen in distinct stages (indexed by $k \in \mathbb{N}$) and guided by the inference results conditioned on data observed from previous experiments. This design setting naturally arises in a range of applications---including medical imaging, subsurface exploration, assessment of environmental hazards, and more---where practitioners need to dynamically adjust experimental conditions, such as sensor locations, based on feedback from previous estimations.
In such settings, sOED's advantage over standard batch experimental design is its ability to further reduce the estimation uncertainty by adaptively tailoring experimental conditions towards the unknown ground truth.
However, sOED tends to be more computationally challenging than batch OED due to the presence of a feedback loop that requires the state of knowledge about the model parameters to be iteratively updated after each experiment is conducted.

To set up the problem, let $\params \in \mathbb{R}^{\Nm}$ denote the unknown model parameters we seek to estimate, $\design_k \in \designSpace_k$ denote the stage-$k$ experimental conditions, and $\data_k \in \dataSpace_k \subset \mathbb{R}^{\Nd}$ denote the corresponding data collected during the $k$-th experiment. \pagebreak
The design $\design_k$ may represent sensor locations (e.g., device placement) or parameterized sources and boundary conditions (e.g., illumination patterns) that drive the system. The design space may be continuous, $\designSpace_k\subset\mathbb{R}^{\Ne}$, or a finite set of $\Ne$ candidate designs. We present our approach in this general setting and defer implementation details to later sections.
We denote the sequence of experimental conditions and the corresponding data in the previous $k{-}1$ experiments by $\Hist_{k{-}1} = [\design_1^*,\data_1^*,\ldots,\design_{k{-}1}^*,\data_{k{-}1}^*]$ (with $\Hist_0 = \emptyset$). At stage-$k$, we represent the posterior density of the inverse problem based on all previous experiments $\Hist_{k{-}1}$ by $\target(\params \given \Hist_{k{-}1})$. This way, after conducting the $k$-th experiment, the posterior density can be recursively updated as
\begin{equation}
\target(\params \given \design_k,\data_k,\Hist_{k{-}1}) = \frac{\likelihood{\data_k}{\params, \design_k, \Hist_{k{-}1}} \,\target(\params \given \Hist_{k{-}1}) }{\target(\data_k \given \design_k, \Hist_{k{-}1})},
\label{eq:posterior_update}
\end{equation}
where $\target(\data_k \given \design_k, \Hist_{k{-}1})$ is the (typically unknown) stage-$k$ evidence at design $\design_k$, and the previous posterior density $\target(\params \given \Hist_{k{-}1})$ becomes the stage-$k$ prior.
Before collecting any data, we initialize the first stage with a prior distribution by setting $\target(\params \given \Hist_0) = \target(\params)$. This initial prior is often chosen to be analytically tractable, for example a Gaussian distribution $\target(\params) \sim \mathcal{N}(\params_0,\Cpr)$. At stage $k$, the likelihood $\likelihood{\data_k}{\params, \design_k, \Hist_{k{-}1}}$ is determined by the model and measurement process that link the parameters and design to the noisy observations.
As an example, we consider an additive Gaussian noise model,
\begin{equation}
\data_k(\design_k) = \PtO(\design_k,\params)+\noise_k(\design_k),
\label{eq:stagek_model}
\end{equation}
where $\PtO \colon \designSpace_k \times \paramSpace \rightarrow \dataSpace_k$ is the forward map, and $\noise_k \sim \mathcal{N}(\bvec{0},\Cnoise(\design_k))$ is the Gaussian noise.
Under these assumptions, the stage-$k$ likelihood is independent of previous experiments, \ie,
\begin{equation}\likelihood{\data_k}{\params, \design_k, \Hist_{k{-}1}} = \likelihood{\data_k}{\params, \design_k} \propto \exp\left(-\frac{1}{2} \norm{\Cnoiseinv{\design_k}}{\PtO(\design_k,\params)-\data_k(\design_k)}^2\right).
\label{eq:likelihood}
\end{equation}
We are particularly interested in problems where the forward map is nonlinear and expensive to evaluate, e.g., those implicitly defined through a partial differential equation (PDE). In these settings, $\params$ typically represents a discretized functional input to the PDE, making it high-dimensional.
For simplicity, we assume constant dimensions for the design and observed data spaces across all experimental stages, and that both the forward map $\PtO$ and noise model remain the same, although these assumptions are not essential to our approach.

There are various sOED formulations, including~\cite{HuanMarzouk:2016:1,ShenHuan:2023:1,FosterIvanovaMalikRainforth:2021:1}, which involve some degree of ``lookahead'', taking future experiments into account when choosing the optimal design at each experimental stage.
We focus on applications without a hard limit on the experimental budget, where a greedy approach to sOED is particularly effective.
At each experimental stage, defining the Kullback--Leibler (KL) divergence of the posterior from the prior
\[
\distKL{\target(\cdot \given \design_k,\data_k,\Hist_{k{-}1})}{\target(\cdot \given \Hist_{k{-}1})} \coloneqq \int \log\left(\frac{\target(\params \given \design_k,\data_k, \Hist_{k{-}1})}{\target(\params \given \Hist_{k{-}1})}\right) \target(\params \given \design_k,\data_k,\Hist_{k{-}1}) \, \mathrm{d}\params,
\]
we choose the design by maximizing the \emph{incremental expected information gain} (iEIG),
\begin{align}
\EIG_k(\design_k) &= \Expect{\data_{k}\given \design_{k},\Hist_{k{-}1}}{\distKL{\target(\cdot \given \design_k,\data_k,\Hist_{k{-}1})}{\target(\cdot \given \Hist_{k{-}1})}}.
\label{eq:EIG_k}
\end{align}
The iEIG~\eqref{eq:EIG_k} does not have a closed-form expression outside specific cases, e.g., those involving linear forward maps, Gaussian priors, and additive Gaussian noise. This leads to a combination of challenges.
The primary challenge arises in approximating nested expectations with respect to the densities $\target(\params \given  \design_k, \data_k, \Hist_{k{-}1})$ and $\target(\data_k \given \design_k,\Hist_{k{-}1})$ in \eqref{eq:EIG_k}, and accessing the ratio $\frac{\target(\params \given \design_k, \data_k, \Hist_{k{-}1})}{\target(\params \given \Hist_{k{-}1})}$ in the KL divergence calculation.
In the first stage, with a tractable prior $\target(\params \given \Hist_{0}) \coloneqq \target(\params)$ that can be directly sampled, the EIG is commonly approximated using a nested Monte Carlo estimator, or a Laplace approximation of the intractable posterior~\cite{BeckMasourEspathTempone:2020,HuanMarzouk:2013:1,WuChenGhattas:2023:1}.
However, approximating the incremental EIGs at subsequent stages ($k > 1$) becomes increasingly complicated, as drawing samples from the intractable prior distribution $\target(\params \given \Hist_{k{-}1})$ in later stages is a known challenge---often requiring specifically designed Markov-chain-based samplers \cite{CuiLawMarzouk:2016:1,girolami2011riemann,martin2012stochastic}. Moreover, the dimensionality of the parameters will further aggravate the complexity of sampling \cite{cotter2013mcmc,mattingly2012diffusion,roberts2001optimal}.
We propose a likelihood-informed, measure-transport-based approach to sequentially approximate optimal designs, which is computationally feasible and scalable with dimension, while simultaneously characterizing the posterior distribution of the inverse problem.

\subsection{Related work}\label{subsec:relatedWork}
The expected information gain (EIG) is most commonly estimated using nested Monte Carlo methods~\cite{Ryan:2003:1,HuanMarzouk:2013:1}, which express the EIG as an expectation of the difference between the log-likelihood and log-evidence. This requires an outer Monte Carlo loop for the expectation and an inner Monte Carlo loop to estimate the intractable evidence for each outer sample.
Reduced-order models have been used to mitigate the computational cost of this nested structure~\cite{WuOLeary-RoseberryChenGhattas:2023:1}. However, while these estimators are asymptotically unbiased, they converge more slowly than standard Monte Carlo and can exhibit substantial finite-sample bias, particularly for concentrated likelihoods~\cite[Section 3.1]{HuanJagalurMarzouk:2024:1}. Recently, measure transport approaches to batch optimal experimental design have been explored in~\cite{CaoChenBrennanOLearyRoseberryMarzoukGhattas:2024:1,LiBaptistaMarzouk:2024:1,DongJacobsenKhalloufiAkramLLiuDuraisamyHuan:2025:1}.
These methods can generally be viewed as two-step estimation approaches, combining transport-based density estimation with Monte Carlo estimation of the expectation.

Our previous work~\cite{KovalHerzogScheichl:2024:1} also uses transport maps for sOED, but is limited to low- to moderate-dimensional parameter spaces. Additional sOED approaches are reviewed in~\cite{HuanJagalurMarzouk:2024:1,RyanDrovandiMcGreePettitt:2016:1}.
Many of these are tailored to low-dimensional parameters and rely on sequential Monte Carlo methods to propagate samples between experiments~\cite{DrovandiMcGreePettitt:2013:1,KleinegesseDrovandiGutmann:2020:1}.
In contrast, we focus on high-dimensional problems where parameters are discretized functional inputs of the forward map.

In the context of OED for large- or infinite-dimensional problems, common approaches typically involve some combination of Gaussian approximations to the posterior~\cite{WuChenGhattas:2023:1}, the use of derivative-informed neural networks~\cite{WuOLeary-RoseberryChenGhattas:2023:1}, or approaches that exploit the presence of low-dimensional structures~\cite{CaoBaptistaChenLiGhattasOdenMarzouk:2023:1,LiBaptistaMarzouk:2024:1}.
The aforementioned works all focus on batch OED, though a combination of these approaches has also been used in sOED.
Specifically,~\cite{GoChen:2024:1} employs Gaussian approximations to the posterior, constructed efficiently using dimension-reduced neural network surrogates, for the sequential selection of optimal observation times in Bayesian inverse problems involving dynamical systems. In contrast, our primary focus is on the optimal selection of sensor locations, which requires a different formulation of the sOED problem.

\subsection{Our approach and contributions}
We propose a novel approach to sequential optimal experimental design by maximizing sharp bounds on the iEIG. Our main result extends the intrinsic dimensionality analysis of Bayesian inverse problems~\cite{LiCuiLiMarzoukZahm:2024:1} to derive a fast-to-evaluate upper bound on the iEIG, which significantly simplifies the computationally demanding sOED objective functions. A key component in constructing this bound in the sequential setting is the use of measure transport. Here we utilize tensor-train techniques to efficiently construct these transport maps, and leverage the likelihood-informed dimension reduction induced by the same iEIG bound to make the transport maps scalable to high- or even infinite-dimensional parameter spaces. Under an additional Gaussian noise assumption, we also provide an alternative bound on the iEIG, and show that both proposed bounds yield tighter estimates of the iEIG than some of the existing EIG surrogates.
Another highlight of our framework is the integration of conditional transport maps into the sOED process to facilitate amortized inference. This enables the construction of transport maps prior to data collection at each experimental stage, allowing for real-time posterior inference once data are available.
Finally, through numerical examples, we evaluate the effectiveness of our designs, comparing them to nested Monte Carlo estimators and assessing the performance of designs maximizing our proposed upper bound against those based on Gaussian approximations.

The rest of this paper is organized as follows. In~\cref{sec:background}, we summarize relevant background material on likelihood-informed parameter dimension reduction, density approximation via triangular measure transport, and functional tensor trains. In~\cref{subsec:iEIG_bounds}, we derive upper bounds on the incremental expected information gain and discuss how these bounds can be used to guide the sequential design of experiments. In~\cref{subsec:sOED_composite}, we discuss how to combine conditional transport maps with likelihood-informed parameter dimension reduction at each experimental stage to enable scalable amortized inference in the sOED procedure.
In~\cref{subsec:sOED_restart}, we improve transport map accuracy by incorporating a restart strategy, leading to our final sOED algorithm, which is then used in numerical examples in~\cref{sec:num1} and~\cref{sec:num2}.

\section{Background}\label{sec:background}

\subsection{Likelihood-informed subspaces}\label{subsec:LIS}
The efficiency of sOED algorithms critically depends on the ability to characterize intractable posterior densities $\target(\params \given \data) \propto \likelihood{\data}{\params}\,\target(\params)$, a task that is particularly challenging in high-dimensions. Most efficient algorithms \cite{BrennanBigoniZahmSpantiniMarzouk:2020:1,ConstantineKentBui-Thanh:2016:1,CuiLawMarzouk:2016:1,CuiMartinMarzoukSolonenSpantini:2014:1,cui2022prior} exploit low-dimensional structure in the parameter–data interaction---arising from prior regularity, smoothness of the forward model, and the incompleteness of the data---to enable scalable posterior inference. We use the \emph{likelihood-informed subspace} (LIS) methods outlined in~\cite{CuiMartinMarzoukSolonenSpantini:2014:1,CuiTong:2022:1,CuiZahm:2021:1,LiCuiLiMarzoukZahm:2024:1,zahm2022certified} to identify such low-dimensional structures.

LIS methods seek to identify a subspace $\mathcal{M}_{\LISdim} \subset \paramSpace$ of dimension $\LISdim \ll \Nm$ that contains the effective support of the likelihood $\likelihood{\data}{\params}$.
Let $\Proj_{\LISdim}$ and $\Proj_{\perp}$ denote orthogonal projection operators onto the space $\paramSubspace_{\LISdim}$ and its orthogonal complement $\paramSubspace_{\perp}$, respectively.
The LIS decomposes the parameter as $\params = \params^{\LISdim} + \params^{\perp}$, with $\params^{\LISdim} = \Proj_{\LISdim}\params$ and $\params^{\perp} = \Proj_{\perp}\params$, and approximates the likelihood by $\likelihood{\data}{\params} \approx \int \likelihood{\data}{\params^{\LISdim}+{\params^{\perp}}'} \target({\params^{\perp}}' \given \params^{\LISdim}) \,\mathrm{d}{\params^{\perp}}' \eqqcolon \Rlikelihood{\data}{\params^{\LISdim}},$
which has support on $\paramSubspace_{\LISdim}$. This results in an approximation to the full-dimensional posterior:%
\begin{equation}
    \target(\params\given \data) \approx \widetilde{\target}(\params \given \data ;\Proj_{\LISdim}) \propto \underbrace{\Rlikelihood{\data}{\params^{\LISdim}} \target(\params^{\LISdim})}_{\widetilde{\target}(\params^{\LISdim} \given \data)}\,{\target(\params^\perp \given \params^{\LISdim})},
    \label{eq:li_posterior}
\end{equation}
where, with a slight abuse of notation, $\target(\params^{\LISdim}) = \int_{\Proj_{\perp}} \target(\params^{\LISdim} + {\params^{\perp}}') \,\mathrm{d} {\params^{\perp}}'$ is the marginal prior and $\target(\params^{\perp}\given\params^{{\LISdim}}) = \frac{\target(\params^{\LISdim}+\params^{\perp})}{\target(\params^{\LISdim})}$ is the conditional prior. Note that $\widetilde{\target}(\params^{\LISdim} \given \data)$ is also the marginal posterior of the reduced parameters. This decomposition offers an efficient two-step sampling strategy: one can sample the lower-dimensional marginal posterior using methods such as Markov chain Monte Carlo, followed by independent samples drawn from the conditional prior. With a slight modification using either the pseudo-marginal principle or importance sampling, exact full posterior samples can also be obtained, see \cite{CuiTong:2022:1,cui2022prior,CuiZahm:2021:1} for details.

Accurate identification of the LIS is central to these algorithms. We outline both \emph{data-dependent} and \emph{data-free} approaches based on the derivative of the likelihood function. Without loss of generality, we assume that the parameters are transformed so that the associated prior distribution is a standard multivariate Gaussian, $\target(\params) = \mathcal{N}(\bvec{0},\I{\Nm})$, see \cite{CuiMartinMarzoukSolonenSpantini:2014:1,martin2012stochastic} for examples of linear transformations and \cite{CuiTong:2022:1} for nonlinear transformations.

\smallskip
\paragraph{Data-dependent LIS~\cite{CuiTong:2022:1}} Given some measured data, $\data \in \dataSpace$, the data-dependent LIS is defined as the subspace spanned by the first $\LISdim$ leading eigenvectors of the Gram matrix,
\begin{equation}
\Hmat(\data) = \int \nabla_{\params} \log \likelihood{\data}{\params} \,\nabla_{\params} \log \likelihood{\data}{\params}^\top \target(\params \given \data) \, \mathrm{d}\params.
\label{eq:gram}
\end{equation}
Denoting the eigenvalues of $\Hmat(\data)$ as $\{\lambda_{i}(\Hmat(\data))\}_{i=1}^{\Nm}$ with $\lambda_1 \geq \lambda_2 \geq \ldots \geq \lambda_{\Nm}$, the approximation to the posterior resulting from a restriction to this data-dependent subspace (with corresponding projection operator $\Proj^{\bvec{y}}_{\LISdim}$) satisfies~\cite[Theorem 2.4]{CuiTong:2022:1}
\begin{equation}
\distH{\target(\cdot \given \data)}{\widetilde{\target}(\cdot \given \data; \Proj^{\bvec{y}}_{\LISdim})} \leq \frac{\sqrt{\kappa}}{2}\bigg( \sum_{k={\LISdim}+1}^{\Nm}\lambda_{k}(\Hmat(\data))\bigg)^{1/2},
\label{eq:DI_bound}
\end{equation}
where $\distH{\cdot}{\cdot}$ is the Hellinger distance and $\kappa$ is the subspace Poincar\'e constant of the prior, which is bounded under mild conditions (see~\cite[Assumption 2.1, Proposition 2.2]{CuiTong:2022:1} for details).

\smallskip
\paragraph{Data-free LIS~\cite{CuiZahm:2021:1,LiCuiLiMarzoukZahm:2024:1}}\label{paragraph:DF-LIS}
Compared to the data-dependent LIS, the data-free LIS is constructed before any data are observed. It identifies directions along which the posterior deviates most from the prior for average data realizations. Given the Fisher information matrix
\begin{equation} \mathcal{I}(\params) = \int \left(\nabla_{\params} \log \likelihood{\data}{\params} \, \nabla_{\params} \log \likelihood{\data}{\params}^\top\right) \likelihood{\data}{\params}\, \mathrm{d}\data,
\label{eq:Fisher}
\end{equation}
the data-free LIS is defined through the dominant $r$-dimensional eigenspace of the matrix
\begin{equation}
\Hmat_{\mathrm{I}} = \int \mathcal{I}(\params) \, \target(\params) \, \mathrm{d}\params.
\label{eq:avgFisher}
\end{equation}
The expected posterior approximation error (over the data distribution) of ~\eqref{eq:li_posterior}, using a projection operator $\Proj^\mathrm{I}_{\LISdim}$ onto the leading eigenspace of $\Hmat_{\mathrm{I}}$, admits a bound analogous to the data-dependent case~\cite{CuiDolgovZahm:2023:1,CuiZahm:2021:1}. With a standard Gaussian prior, this bound can be further sharpened~\cite[Theorem~9]{LiCuiLiMarzoukZahm:2024:1} to
\begin{equation}
    \mathbb{E}_{\data}\left[ \distKL{\target(\cdot\given \data)}{\widetilde{\target}(\cdot \given \data ; \Proj^\mathrm{I}_{\LISdim})}\right] \leq \frac{1}{2}\sum_{k=\LISdim+1}^{\Nm} \log \left(1+\lambda_k(\Hmat_{\mathrm{I}})\right).
    \label{eq:EIG_bound}
\end{equation}

To apply these approaches, one needs to numerically approximate the matrices $\Hmat(\bvec{y})$ and $\Hmat_{\mathrm{I}}$, which typically involves Monte Carlo approximations to the expectation. As shown in~\cite{CuiZahm:2021:1}, the data-dependent approach provides more accurate approximations to the informed subspace for a given instance of data. However, it is computationally challenging due to the calculation of the expectations over the typically inaccessible posterior. The data-free approach avoids this challenge and is advantageous when solving multiple Bayesian inverse problems with varying observed data. Furthermore, the average Fisher information $\Hmat_{\mathrm{I}}$ in~\eqref{eq:EIG_bound} immediately provides a sharp bound on the expected information gain, offering a fast surrogate to guide our sOED. Specifically, setting $r=0$ in~\eqref{eq:EIG_bound} gives
\begin{equation}
    \mathbb{E}_{\data}\left[ \distKL{\target(\cdot\given \data)}{{\target}(\cdot )}\right] \leq \frac{1}{2} \log \det \left(\I{\Nm}+\Hmat_{\mathrm{I}}\right).
    \label{eq:EIG_bound_full}
\end{equation}

\subsection{Knothe--Rosenblatt rearrangement}
Although the LIS methods can reduce the dimensionality of the parameters, the question of how to efficiently characterize and sample from the dimension-reduced posterior remains. Measure transport, see~\cite{BaptistaMarzoukZahm:2023:1,MoselhyMarzouk:2012:1,MarzoukMoselhyParnoSpantini:2016:1} for instance, offers a versatile solution to this, with applications widely found in Bayesian inference, rare event estimation, and optimal experimental design. Measure transport constructs an invertible transformation $\mathcal{S} \colon \mathbb{R}^{n} \rightarrow \mathbb{R}^{n}$ between a tractable reference measure $\nu$ with density $\refD$ (\eg, a multivariate Gaussian) and the intractable target measure $\mu$ with density $\target$ (\eg, a posterior density).
If $\mathbb{R}^{n} \ni \refparams = \mathcal{S}(\targetparams) \sim \refD$ for any $\mathbb{R}^n \ni \targetparams \sim \target$ then $\mathcal{S}$ is said to \emph{push forward} the target density $\target$ to the reference density $\refD$ (likewise, \emph{pull back} $\refD$ to $\target$). The pushforward and pullback operators are defined, respectively, as
\begin{align*}
\pushforward{\mathcal{S}}\target(\refparams) &= \left(\target \circ \mathcal{S}^{-1}\right)(\refparams) \det\big( \Jac\mathcal{S}^{-1} (\refparams)\big) = \refD(\refparams) \\
\pullback{\mathcal{S}}\refD(\targetparams) &= \left(\refD \circ \mathcal{S}\right)(\targetparams) \det\big( \Jac\mathcal{S} (\targetparams) \big) = \target(\targetparams),
\end{align*}
where $\Jac\mathcal{S}$ denotes the Jacobian of $\mathcal{S}$.
In our setting, we assume that both the reference and target measures are absolutely continuous with respect to the Lebesgue measure.

For multivariate random variables, there may exist infinitely many maps $\mathcal{S}$ that couple the reference and the target distributions and various numerical implementations, e.g., ~\cite{chen2018neural,KruseDetommasoKoetheScheichl:2021:1,onken2021ot,PapamakariosNalisnickRezendeMohamedLakshminarayanan:2021:1}. Our sOED approach is flexible to utilize various transport maps. In this work, we employ the \emph{Knothe-Rosenblatt} (KR) rearrangement:
\begin{equation}
	\refparams
	=
	\begin{bmatrix}
		\refparam_1
		\\
		\refparam_2
		\\
		\vdots
		\\
		\refparam_n
	\end{bmatrix}
	=
    \mathcal{S}(\targetparams)
    =
	\begin{bmatrix*}[l]
		\mathcal{S}_{\targetparam_1}(\targetparam_1)
		\\
		\mathcal{S}_{\targetparam_2 \given \targetparam_1}(\targetparams_{1:2})
		\\
		\quad
		\vdots
		\\
		\mathcal{S}_{\targetparam_n \given \targetparams_{1:n-1}}(\targetparams)
	\end{bmatrix*}
	,
	\label{eq:KR_map_x}
\end{equation}
where $\targetparams_{1:k} = [\targetparam_1,\ldots,\targetparam_k]^\top$.
Each component of the KR map, $\mathcal{S}_{\targetparam_k \given \targetparams_{1:k{-}1}}(\targetparams_{1:k{-}1},\targetparam_k)$, is monotone in the last input parameter $\targetparam_k$. 
Given a sample from the reference density, $\refparams \sim \refD$, a sample from the target density, $\targetparams \sim \target$, can be obtained via the inverse map $\cT \coloneqq \mathcal{S}^{-1}$.
The triangular structure of the KR map enables application of $\cT$ via sequential inversion of univariate functions. 
Additionally, the components of the KR map are defined through the marginal conditionals of the target density. Thus, the $k$-th component of $\mathcal{S}$ immediately grants access to the conditional target density $\target(\targetparam_k \given \targetparams_{1:k{-}1})$.
This feature makes the KR map particularly valuable for Bayesian inference and, as we will detail in~\cref{sec:sOED}, for sOED.

\paragraph{Tensor-train construction} \label{subsec:DIRT} We use the density approximations strategy of~\cite{CuiDolgov:2022:1,CuiDolgovZahm:2023:2,dolgov2020approximation,WestermannZech:2023;1}, particularly the squared tensor trains~\cite{CuiDolgov:2022:1}, to implement KR maps. Consider a target density $\target = \frac{\targetUnnorm}{Z}$, where $\targetUnnorm$ can be evaluated pointwise and the normalizing constant $Z$ is unknown. We approximate the square root of the unnormalized density $\targetUnnorm$ using a functional tensor train,
\begin{equation}
    \sqrt{\targetUnnorm(\targetparams)} \approx \TTCore{1}{\targetparam_1} \cdots \TTCore{i}{\targetparam_i} \cdots \TTCore{n}{\targetparam_n}  \eqqcolon \widehat{f},
    \label{eq:FTT}
\end{equation}
where each $\TTCore{i}{\targetparam_i} \in \mathbb{R}^{{\TTRank_{i-1}}\times{\TTRank_i}}$ is a matrix-valued univariable function, with $\TTRank_{0} = \TTRank_{n} = 1$.
The elements of each matrix-valued function, $[\TTCore{i}{\targetparam_i}]_{k,j}$ are represented as a linear combination of $M_i$ basis functions.
Such tensor train approximations can be constructed efficiently using alternating-direction cross interpolation methods~\cite{GorodetskyKaramanMarzouk:2019:1,OseledetsTyrtyshnikov:2010:1}.
Here we employ a functional extension of the rank-adaptive alternating minimal energy scheme as described in~\cite[Appendix B]{CuiDolgov:2022:1} to build tensor trains.
Denoting $M {=} \max_{i} M_i$ and $r {=} \max_{i} \TTRank_i$, it requires $\mathcal{O}(nM\TTRank^2)$ density evaluations to construct such approximations.
The tensor-train approximation $\widehat{f}$ leads to a normalized approximated target density
\begin{equation}
    \widehat{\target}(\targetparams) = \frac{\widehat{\targetUnnorm}(\targetparams)}{\xi}, \quad \widehat{\targetUnnorm}(\targetparams) = \widehat{f}(\targetparams)^2+\tau\refD(\targetparams), \quad \xi = \int \widehat{\targetUnnorm}(\targetparams) \, \mathrm{d}\targetparams,
\end{equation}
where $\refD(\targetparams) = \prod_{i=1}^n \refD_{i}(\targetparam_i)$ is a product-form reference density such that $\sup_\targetparams \frac{\targetUnnorm(\targetparams)}{\refD(\targetparams)} < \infty$ and $\tau > 0$ is a constant. The ``defensive'' term $\tau \refD$ is added to ensure that the tensor train surrogate $\widehat{\target}$ can define importance sampling estimators satisfying the central limit theorem. As shown in~\cite[Theorem 1]{CuiDolgov:2022:1}, choosing tensor ranks $\TTRank_i$ to bound the approximation error by some threshold $\varepsilon$, i.e., $\|\widehat{f} - \sqrt{\targetUnnorm}\| < \varepsilon$, and using a constant $\tau < \varepsilon$ ensures that the Hellinger error of the approximate posterior is bounded such that $\distH{\target}{\widehat{\target}} \leq \frac{2\varepsilon}{\sqrt{Z}}$. This, in turn, bounds errors in expectations computed over the approximate posterior.

Leveraging the separable structure of tensor trains, the marginal densities
\begin{equation}
\widehat{\target}_{\targetparams_{1:k}}(\targetparams_{1:k}) = \int_{\targetSpace_{k+1:n}} \widehat{\target}(\targetparams) \,\mathrm{d}\targetparams_{k+1:n}, \quad k = 1,\ldots,n-1,
\end{equation}
of the approximation $\widehat{\target}$ can be computed dimension-by-dimension with $\mathcal{O}(nM\TTRank^3)$ floating point operations. This naturally defines a sequence of distribution functions
\[
\mathcal{S}_{\targetparam_k\given\targetparams_{1:k{-}1}}(\targetparam_k) = \int_{-\infty}^{x_k} \widehat{\target}_{\targetparam_k \given \targetparams_{1:k{-}1}} \, \mathrm{d}\targetparam_k' = \int_{-\infty}^{x_k} \frac{\widehat{\target}_{\targetparams_{1:k}}(\targetparams_{1:k{-}1},\targetparam_k')}{\widehat{\target}_{\targetparams_{1:k{-}1}}(\targetparams_{1:k{-}1})} \, \mathrm{d}\targetparam_k', \quad k = 2, \ldots, n,
\]
together with $\mathcal{S}_{\targetparam_1}(\targetparam_1) {=} \int_{-\infty}^{\targetparam_1} \widehat{\target}_{\targetparam_1}(\targetparam_1')\,\mathrm{d}{\targetparam_1'}$, this gives a KR map $\mathcal{S}$ coupling the approximation $\widehat{\target}$ with the uniform density $\refD_{\text{uni}}$ on $(0,1]^n$, \ie, $\pushforward{\mathcal{S}}\widehat{\target} = \refD_{\text{uni}}$.
Since the product-form reference $\refD$ admits a diagonal KR rearrangement, $\mathcal{R} = [\mathcal{R}_{1}(\refparam_1),\ldots,\mathcal{R}_n(\refparam_n)]^\top$ with $\mathcal{R}_{i}(\refparam_{i}) {=} \int_{-\infty}^{\refparam_i} \refD_{i}(\refparam_i')\,\mathrm{d}\refparam_i'$, such that $\pushforward{\mathcal{R}}\refD = \refD_{\text{unif}}$, the composite map $\cT = \mathcal{S}^{-1} \circ \mathcal{R}$ satisfies $\pushforward{\cT}\refD = \widehat{\target}$. 

\smallskip
\paragraph{Deep approximation} 
For targets that are highly concentrated or exhibit strongly nonlinear dependence, accurate tensor-train approximations may require large ranks and many basis functions. Since the number of unnormalized density evaluations scales as $\mathcal{O}(nM\TTRank^2)$, the efficiency depends critically on the ranks and the number of basis functions, especially when each density evaluation requires a PDE solve.
The deep approximation approach~\cite{CuiDolgov:2022:1} addresses this challenge by constructing a composite map
\(
\cT^{L}=\cQ^{1} {\circ} \cQ^{2} {\circ} \cdots {\circ} \cQ^{L},
\)
where each $\cQ^\ell$ is easier to compute. The composition is built recursively using a sequence of bridging densities $\{\target^{\ell}\}_{\ell=1}^{L}$, with $\target^{L}=\target$, that progressively capture the complexity of the target.
At layer $\ell$, suppose the previous composition $\cT^{\ell-1}=\cQ^{1}{\circ}\cdots{\circ}\cQ^{\ell-1}$ satisfies $(\cT^{\ell-1})_\sharp \refD \approx \target^{\ell-1}$. Then the pullback of the next bridging density $\target^{\ell}$ under $\cT^{\ell-1}$ yields an approximation
\[
\pullback{(\cT^{\ell{-}1})}\target^\ell \approx \frac{\target^{\ell} \circ \cT^{\ell{-}1}}{\target^{\ell{-}1} {\circ} \cT^{\ell{-}1}} \, \refD,
\]
which is the perturbation of the reference density $\refD$ by the ratio $\frac{\target^{\ell} \circ \cT^{\ell{-}1}}{\target^{\ell{-}1} \circ \cT^{\ell{-}1}}$.
For appropriately chosen bridging densities, the pullback density $\pullback{(\cT^{\ell{-}1})}\target^\ell$ is easier to approximate than the target density itself using tensor trains with lower ranks. We then build the intermediate map $\cQ^{\ell}$ to couple the reference density with the pullback density, $\pushforward{(\cQ^{\ell})}\refD \approx \pullback{(\cT^{\ell{-}1})}\target^\ell$ and enrich the composition as $\cT^{\ell} = \cT^{\ell{-}1}\circ\cQ^{\ell}$. We refer the readers to~\cite{CuiDolgov:2022:1} for further details.

\section{Fast sOED and scalable amortized inference}\label{sec:sOED}
We first extend the upper bound~\eqref{eq:EIG_bound_full} to the iEIG setting to guide sequential design, and show that it also induces dimension reduction for the selected experimental condition. We also introduce an alternative iEIG bound under additional assumptions and discuss connections among various iEIG surrogates. We then combine conditional KR maps with the resulting dimension reduction to enable scalable amortized posterior inference, and propose a restart procedure to improve accuracy.

\subsection{Fast sOED via iEIG bounds}\label{subsec:iEIG_bounds}
\subsubsection{Setup}\label{subsec:iEIG_setup}
At stage $k$, sOED aims to maximize the iEIG~\eqref{eq:EIG_k} from the prior---which is the previous posterior $\target(\params \given \Hist_{k{-}1})$---to the updated posterior $\target(\params \given \design_k, \data_k,\Hist_{k{-}1})$. We consider there exists a measure transport ${\cT^{k{-}1}} {\colon} \mathbb{R}^{\Nm} {\rightarrow} \mathbb{R}^{\Nm}$ that maps the standard Gaussian reference $\refD(\refparams)$ to the stage-$k$ prior $\target(\params \given \Hist_{k{-}1})$, that is $\pushforward{({{\cT^{k{-}1}}})} \refD(\params) = \target(\params \given \Hist_{k{-}1})$.
We define the composite forward map
\begin{equation}
\mathcal{J}(\design_k,\refparams) =  \PtO(\design_k,\cT^{k{-}1} (\refparams) ),\label{eq:compose_forward}
\end{equation}
in the stage-$k$ experiment, and denote the Jacobian of the composite forward map by
\[
\Jac_\refparams \mathcal{J}(\design_k,\refparams) = \Jac_\params \PtO(\design_k,\params)\,\Jac_\refparams{{\cT^{k{-}1}}(\refparams)} , \quad \params = {\cT^{k{-}1}}(\refparams),
\]
where $\Jac_\params \PtO(\design_k,\params) \in \mathbb{R}^{\Nd \times \Nm}$ denotes the Jacobian of the forward map evaluated at $(\design_k,\params)$.

\smallskip
\paragraph{Candidate design setting}
Although the iEIG surrogates discussed below apply to both continuous and discrete design spaces, we focus on a discrete candidate set $\designSpace_k = \{\bvec{E}_1,\ldots,\bvec{E}_{\Ne}\}$ with cardinality $\Ne$, where the goal is to select the optimal experimental condition from this set.
In this context, it is convenient to separate the design from the forward map by introducing the ``stacked'' forward map $\mathcal{F}\colon\mathbb{R}^{\Nm} \rightarrow \mathbb{R}^{\Ne\Nd}$,
\begin{equation}
\mathcal{F}(\params) = [\mathcal{G}(\bvec{E}_1,\params);\mathcal{G}(\bvec{E}_2,\params);\ldots;\mathcal{G}(\bvec{E}_{\Ne},\params)]^\top,
\label{eq:fullPtO}
\end{equation}
that collects the observables over the entire candidate design set. The observables associated with a specific design $\design_k \in \designSpace_k$ can be extracted using a row selection matrix $\bvec{W}(\design_k)$.
Specifically, let $\design_k = \bvec{E}_{j}$, then $\bvec{W}(\design_k) = \unitvec^{\Ne}_j \otimes \I{\Nd}$ where $\unitvec^{\Ne}_j \in \{0,1\}^{1\times\Ne}$ is a row vector corresponding to the $j$-th canonical basis, i.e., its only non-zero entry is the $j$-th element, $\I{\Nd}$ is the $\Nd$-dimensional identity matrix, and $\otimes$ is the Kronecker product.
The design can be separated from the forward map as
\(
\PtO(\design_k,\params) = \bvec{W}(\design_k)\,\mathcal{F}(\params).
\)
Consequently, we have
\[
\mathcal{J}(\design_k, \refparams) = \bvec{W}(\design_k)\, \mathcal{F}\circ \cT^{k{-}1}(\refparams) \quad \text{and} \quad \Jac_\refparams \mathcal{J}(\design_k,\refparams) = \bvec{W}(\design_k)\,\Jac_\refparams ( \mathcal{F} \circ \cT^{k{-}1}(\refparams) ),
\]
where $\Jac_\refparams ( \mathcal{F} \circ \cT^{k{-}1}(\refparams) ) \in \mathbb{R}^{\Ne\Nd \times \Nm}$.

\subsubsection{An information-based bound}\label{subsec:iEIG_UB}
As a start, we consider the pullback density of the stage-$k$ posterior \eqref{eq:posterior_update} under the transport map $\cT^{k{-}1}$,
\begin{align}
\ratio{k}(\refparams \given \design_k,\data_k) = \pullback{({\cT^{k{-}1}})}\target(\refparams \given \design_k, \data_k,\Hist_{k{-}1}) &\propto {\likelihood{\data_k}{{\cT^{k{-}1}}(\refparams),\design_k}\,\refD(\refparams)}.
\label{eq:ratioFun}
\end{align}
which can be viewed as the posterior arising from solving the Bayesian inverse problem with a likelihood $\likelihood{\data_k}{{\cT^{k{-}1}}(\refparams),\design_k}$ and a standard multivariate Gaussian prior. The corresponding Fisher information matrix $\mathcal{I}(\design_k,\cT^{k{-}1},\refparams)$ at a transformed parameters $\refparams$ becomes
\begin{align}
    \mathcal{I}(\design_k,\cT^{k{-}1}\!,\refparams)
    \! &= \hspace{-4pt} \int \!\! \left( \nabla_{\refparams} \log \likelihood{\data_k}{{\cT^{k{-}1}}(\refparams),\design_k} \,\nabla_{\refparams} \log \likelihood{\data_k}{{\cT^{k{-}1}}(\refparams),\design_k}^\top
    \! \right) \! \likelihood{\data}{\cT^{k{-}1}(\refparams),\design_k} \mathrm{d}\data \nonumber \\
    \! & = \Jac{{\cT^{k{-}1}}(\refparams)}^\top \mathcal{I}(\design_k,\params) \, \Jac{{\cT^{k{-}1}}(\refparams)},
    \label{eq:avg_iFisher}
\end{align}
where
$
\mathcal{I}(\design_k,\params) = \Expect{\data \given \params,\design_k}{\nabla_{\params} \log \likelihood{\data_k}{\bvec{m},\design_k} \, \nabla_{\params} \log \likelihood{\data_k}{\params,\design_k}^\top }
$
with $\params = {\cT^{k{-}1}}(\refparams)$ is the Fisher information at the original parameter $\params$.

Since the KL divergence is invariant to invertible transforms, we can then apply the transport map $\cT^{k{-}1}$ to rewrite the iEIG~\eqref{eq:EIG_k} as
\begin{align}
\EIG_k(\design_k)
&= \Expect{\data_{k}\given \design_{k},\Hist_{k{-}1}}{\distKL{\pullback{({{\cT^{k{-}1}}})}\target(\cdot \given \design_k, \data_k,\Hist_{k{-}1})}{\pullback{({{\cT^{k{-}1}}})}\target(\cdot \given \Hist_{k{-}1})}} \nonumber \\
&= \Expect{\data_{k}\given \design_{k},\Hist_{k{-}1}}{\distKL{\ratio{k}(\cdot \given \design_k, \data_k)}{\refD(\cdot)}}.
\label{eq:iEIG_pullback}
\end{align}
Thus, we can apply the bound~\eqref{eq:EIG_bound_full} to the transformed posterior $\ratio{k}$ to bound the iEIG as
\begin{equation}
\EIG_k(\design_k) \leq \EIG_k^\mathrm{I}(\design_k) := \frac{1}{2} \log \det \big( \I{\Nm} + \Hmat_{\mathrm{I}}(\design_k,\cT^{k{-}1}) \big) , \label{eq:iEIG_uB}
\end{equation}
where $\Hmat_{\mathrm{I}}(\design_k,\cT^{k{-}1}) = \int \mathcal{I}(\design_k,\cT^{k{-}1},\refparams) \, \refD(\refparams)\,\mathrm{d}\refparams \in \mathbb{R}^{\Nm \times \Nm}$. This upper bound can serve as a guide for approximating sequentially optimal designs.

For many likelihood functions, low-rank decompositions of the Fisher information matrix $\mathcal{I}(\design_k,\cT^{k{-}1},\refparams)$ can be constructed at each $\refparams$ by combining iterative and/or randomized linear algebra algorithms and adjoint-based evaluations of matrix–vector products. Using our example of additive Gaussian likelihood~\eqref{eq:likelihood} and the composite forward map $\mathcal{J}(\design_k,\refparams) =  \PtO(\design_k,\cT^{k{-}1} (\refparams) )$, the Fisher information matrix takes the form
\begin{align}
     \mathcal{I}(\design_k,\cT^{k{-}1},\refparams) = \Jac_\refparams \mathcal{J}(\design_k,\refparams)^\top \Cnoiseinv{\design_k}\,\Jac_\refparams \mathcal{J}(\design_k,\refparams) ,
     \label{eq:Fisher_gauss}
\end{align}
where the singular value decomposition of the Jacobian $\Jac_{\refparams}\mathcal{J}$ can be computed using randomized methods~\cite{martinsson2020randomized}. The number of adjoint PDE solves, necessary  to compute this to sufficient accuracy depends primarily on the singular value spectrum of the Jacobian. This may scale well when the candidate design set is large.
This way, we can then approximate the stage-$k$ average information matrix $\Hmat_{\mathrm{I}}(\design_k,\cT^{k{-}1})$ using Monte Carlo integration,
\begin{equation}
    \overline{\Hmat}_{\mathrm{I}}(\design_k,\cT^{k{-}1}) = \frac{1}{N} \sum_{i=1}^N  \mathcal{I}(\design_k,\cT^{k{-}1},\refparams^{(i)}),
    \label{eq:MC_Fisher}
\end{equation}
which is subsequently used to find the approximate optimal design
\begin{equation}
\design_k^* \in \Argmax_{\design_k} \overline{\EIG}_k^{\mathrm{I}}(\design_k) := \Argmax_{\design_k} \frac{1}{2}\log \det \big( \I{\Nm} +  \overline{\Hmat}_{\mathrm{I}}(\design_k,\cT^{k{-}1}) \big)  .
\label{eq:iUB}
\end{equation}

\paragraph{Candidate design setting}\label{paragraph:candidate_designs}

The estimator $\overline{\Hmat}_{\mathrm{I}}(\design_k,\cT^{k{-}1})$ may need to be re-evaluated at each new design $\design_k$ in maximizing $\overline{\EIG}_k^{\mathrm{I}}(\design_k)$. This cost can be reduced in the candidate design setting. In particular, using the stacked forward map \eqref{eq:fullPtO}, the estimator $\overline{\Hmat}_{\mathrm{I}}(\design_k,\cT^{k{-}1})$ becomes
\begin{equation}
\overline{\Hmat}_{\mathrm{I}}(\design_k,\cT^{k{-}1}) \!=\! \frac{1}{N} \!\sum_{i=1}^N \! \Jac_\refparams ( \mathcal{F}\circ \cT^{k{-}1}(\refparams^{(i)}))^\top \bvec{W}(\design_k)^\top \Cnoiseinv{\design_k}\bvec{W}(\design_k) \,\Jac_\refparams ( \mathcal{F}\circ \cT^{k{-}1}(\refparams^{(i)})).
    \label{eq:MC_Fisher_design}
\end{equation}
We can pre-compute the Jacobian $\Jac_\refparams ( \mathcal{F}\circ \cT^{k{-}1}(\refparams^{(i)}))$ at a set of $N$ samples $\refparams^{(i)} \sim \refD$ at the start of each design stage. Then, the function $\overline{\EIG}_k^{\mathrm{I}}(\design_k)$ can be evaluated at each candidate design $\design_k \in \designSpace_k$ by extracting the relevant rows of the Jacobian using the selection matrix $\bvec{W}(\design_k)$, without additional computation.
As a by-product, the eigendecomposition of the average information matrix at the optimal design, $\overline{\Hmat}_{\mathrm{I}}(\design^*,\cT^{k{-}1})$, naturally results in the data-free LIS, as detailed in~\cref{subsec:LIS}. The procedure is outlined in~\cref{Alg:DLIS_combined}.
In later sections, the resulting LIS is used to construct scalable conditional maps, enabling amortized inference and providing the required transport maps for sOED in the new stage.

\subsubsection{A covariance-based bound for additive Gaussian noise}\label{subsec:iEIG_EB}

Recall that for a random variable $\bvec{x} \in \mathbb{R}^{n}$ with density $\pi(\bvec{x})$ and bounded second  moment, its differential entropy $h(\bvec{x}) = -\Expect{}{\log \pi(\bvec{x})}$ has the upper bound~\cite[Theorem 9.6.5]{Cover:1999:1}
\begin{equation}
h(\bvec{x}) \leq \frac{1}{2}\Big(n \log (2\pi e) + \log \det(\Cov(\bvec{x})) \Big),
\label{eq:entropy_bound}
\end{equation}
with equality attained when $\bvec{x}$ is Gaussian. We use this result to derive an alternative, gradient-free bound on the iEIG in the case of additive Gaussian observation noise.

At the $k$-th experimental stage, the entropy of the marginal data $\data_k \given \design_k, \Hist_{k{-}1}$ is given by
\[
h(\data_k \given \design_k, \Hist_{k{-}1}) = -\Expect{\data_k \given \design_k, \Hist_{k{-}1}}{\log\target(\data_k \given \design_k, \Hist_{k{-}1})},
\]
and the entropy of the conditional data $\data_k \given \design_k, \Hist_{k{-}1}$ is given by
\[
h(\data_k \given \params, \design_k,\Hist_{k{-}1}) = - \Expect{\data_k \given \params, \design_k, \Hist_{k{-}1}}{\log \likelihood{\data_k}{\params,\design_k,\Hist_{k{-}1}}}.
\]
Following these definitions, the iEIG~\eqref{eq:EIG_k} can also be written as the difference between the entropy of $\data_k \given \design_k, \Hist_{k{-}1}$ and the expected entropy of $\data_k \given \params, \design_k, \Hist_{k{-}1}$; that is,
\begin{align}
\EIG_k(\design_k) &= h(\data_k \given \design_k, \Hist_{k{-}1})-\Expect{\params \given \Hist_{k{-}1}}{h(\data_k \given \params, \design_k,\Hist_{k{-}1})}.\label{eq:iEIG_entropy}
\end{align}
For additive Gaussian noise, the likelihood is Gaussian with covariance matrix $\Cnoise(\design_k)$, thus the entropy of $\data_k \given \params, \design_k, \Hist_{k{-}1}$ has the closed form expression, $$
h(\data_k \given \params, \design_k, \Hist_{k{-}1})= \frac{1}{2} \Big( \Nd \log (2\pi e) + \log \det(\Cnoise(\design_k)) \Big).
$$
which is independent of the parameter $\params$. Plugging this into \eqref{eq:iEIG_entropy}, and making use of~\eqref{eq:entropy_bound} to bound the entropy $h(\data_k \given \design_k, \Hist_{k{-}1})$, we have the following covariance bound on the iEIG,
\begin{align*}
\EIG_k(\design_k) &\leq \frac{1}{2}\log\det(\Cov_{\data_k\given\design_k, \Hist_{k{-}1}}(\data_k)) - \frac{1}{2}\log\det(\Cnoise(\design_k)) \eqqcolon \EIG_k^{\mathrm{C}}(\design_k),
\end{align*}
where \(
\Cov_{\data_k\given\design_k, \Hist_{k{-}1}}(\data_k) = \Cov_{\params \given \bvec{H}_{k{-}1} }\big(\PtO(\design_k,\params)\big) + \Cnoise(\design_k)
\)
as the new data are independent of previous experiments. To allow for straightforward numerical evaluation, we apply the transport map $\cT^{k{-}1}$ to rewrite $\Cov_{\params \given \bvec{H}_{k{-}1}}(\PtO(\design_k,\params) )$ as the covariance of the composite forward map $\mathcal{J}(\design_k,\refparams) =  \PtO(\design_k,\cT^{k{-}1} (\refparams) )$.
This simplifies the bound to
\begin{align}
\EIG_k^{\mathrm{C}}(\design_k) & = \frac{1}{2}\log\det\big(\Cov_{\bvec{z}}\big( \mathcal{J}(\design_k,\refparams)\big) + \Cnoise(\design_k) \big) - \frac{1}{2}\log\det(\Cnoise(\design_k))\nonumber \\
& = \frac{1}{2}\log\det\big(\I{\Nd} + \Cnoise(\design_k)^{-1} \Cov_{\bvec{z}}\big( \mathcal{J}(\design_k,\refparams)\big) \big). \label{eq:entropic_bound}
\end{align}

\paragraph{Candidate design setting} In this case, the covariance bound simplifies to
\[
\EIG_k^{\mathrm{C}}(\design_k) = \frac{1}{2}\log\det\big(\I{\Nd} + \Cnoise(\design_k)^{-1} \bvec{W}(\design_k)\Cov_{\refparams}\big(\mathcal{F}\circ \cT^{k{-}1}(\refparams)\big) \bvec{W}(\design_k)^{\top} \big),
\]
where the design enters through the selection matrix $\bvec{W}(\design_k)$.
This permits us to pre-compute $\Cov_{\refparams}(\mathcal{F}\circ \cT^{k{-}1}(\refparams))$ at the start of the stage-$k$ experiment.
Using reference samples $\refparams^{(i)} \sim \refD(\cdot)$ for $i = 1, \ldots, N$, the covariance $\Cov_{\refparams}(\mathcal{F}\circ \cT^{k{-}1}(\refparams))$ has the Monte Carlo estimator
\begin{equation}
\overline{\bvecS{\Gamma}}_k(\mathcal{F}) \!:=\! \frac{1}{N{-}1} \! \sum_{i=1}^N \! \left( \mathcal{F}{\circ} \cT^{k{-}1}(\refparams^{(i)}) {-} \bar{\mu}\right) \!\! \left( \mathcal{F}{\circ} \cT^{k{-}1}(\refparams^{(i)}) {-} \bar{\mu}\right)^{\!\!\top}\!\!\!, \quad \bar{\mu} \!=\! \frac{1}{N} \! \sum_{i=1}^N \! \mathcal{F}{\circ} \cT^{k{-}1}(\refparams^{(i)}),
\label{eq:Cov_estimator}
\end{equation}
which leads to the approximation of the optimal design via the covariance bound,
\begin{equation}
\design_k^* \in \Argmax_{\design_k} \overline{\EIG}_k^{\mathrm{C}}(\design_k) := \Argmax_{\design_k} \frac{1}{2}\log\det\big(\I{\Nd} + \Cnoise(\design_k)^{-1} \bvec{W}(\design_k)\overline{\bvecS{\Gamma}}_k(\mathcal{F})\bvec{W}(\design_k)^{\top}\big).
\label{eq:iUB_E}
\end{equation}
Unlike the information-based bound, this bound cannot be used to construct an LIS. To develop scalable algorithms for building transport maps, we therefore estimate the LIS using the averaged information matrix $\overline{\Hmat}_{\mathrm{I}}(\design_k^*,\cT^{k-1})$ evaluated at the optimal design. The procedure for optimizing the covariance bound~\eqref{eq:iUB_E} is summarized in~\cref{Alg:DLIS_combined}.

\begin{algorithm}[h]
\caption{Sensor selection using iEIG surrogates and construction of the data-free LIS basis $\LISDF$ at the optimal design $\design^*$. The inputs are a transport map $\cT$, a stacked forward map $\mathcal{F}$, a truncation tolerance $\epsilon_I$, and a sample size $N$. \label{Alg:DLIS_combined}}
\begin{algorithmic}[1]
\Function{iEIGUB}{$\cT,\mathcal{F},\epsilon_I,N$}
\If{using information-based surrogate~\eqref{eq:iUB}}
\State Draw samples $\refparams^{(i)} \sim \refD$ and compute $\Jac_\refparams ( \mathcal{F}\circ \cT(\refparams^{(i)}))$ for $i=1,\ldots,N$
\State $\design^* \gets \Argmax_{\design \in \designSpace} \frac{1}{2}\log\det(\I{\Nm}+\overline{\Hmat}_{\mathrm{I}}(\design,\cT))$, where $\overline{\Hmat}_{\mathrm{I}}(\design,\cT)$ is given in~\eqref{eq:MC_Fisher_design}
\EndIf
\If{using covariance-based surrogate~\eqref{eq:iUB_E}}
\State Draw samples $\refparams^{(i)} \sim \refD$ and compute $\mathcal{F}\circ \cT(\refparams^{(i)})$ for $i=1,\ldots,N$
\State $\design^* \gets \Argmax_{\design \in \designSpace} \frac{1}{2}\log\det(\I{\Nd}+ \Cnoise(\design)^{-1} \bvec{W}(\design) \overline{\bvecS{\Gamma}}_k(\mathcal{F}) \bvec{W}(\design)^{\top})$ using~\eqref{eq:Cov_estimator}
\State Build $\overline{\Hmat}_{\mathrm{I}}(\design^*,\cT)$ as in~\eqref{eq:MC_Fisher_design}
\EndIf
\State Compute eigendecomposition $\overline{\Hmat}_{\mathrm{I}}(\design^*,\cT) = \LISDF \bvecS{\Lambda} \LISDF^\top$
\State Set $\LISDF = \LISDF(:,1:\LISdim)$, with $\LISdim$ chosen to ensure $\frac{1}{2}\sum_{i=\LISdim+1}^{\Nm}\log(1+\lambda_i(\overline{\Hmat}_\mathrm{I}(\design^*,\cT))) \leq \mathrm{\epsilon_I}$
\State \Return $\design^*,\LISDF$
\EndFunction
\end{algorithmic}
\end{algorithm}

\subsubsection{Remarks on iEIG surrogates} Let $\mathbb{S}^n_{+}$ denote the cone of symmetric positive semi-definite matrices and let $\bvecS{\Gamma} \in \mathbb{S}^n_{+}$ be positive definite. The following identities are used:
\begin{itemize}[leftmargin=19pt]
\item[(\romannumeral 1)] The map
\(
\varphi(\bvec{A}) \coloneqq \log\det\!\bigl(\bvec{I}+\bvecS{\Gamma}^{-1}\bvec{A}\bigr)
\)
is order-preserving on $\mathbb{S}^n_{+}$  with respect to the Loewner partial order~\cite[Lemma~2.3]{tie2011rearrangement}, i.e., $\varphi(\bvec{A})\leq\varphi(\bvec{B})$ for all $\bvec{A}\preceq \bvec{B}$.
\item[(\romannumeral 2)] The difference
\(
g(\bvec{A}) = \mathrm{trace} (\bvecS{\Gamma}^{-1} \bvec{A}) - \log\det( \bvec{I} + \bvecS{\Gamma}^{-1} \bvec{A} )
\)
satisfies $g(\bvec{A}) \geq 0$ for all $\bvec{A} \in \mathbb{S}^n_{+}$, with strict inequality whenever $\bvec{A}\neq \bvec{0}$. Moreover, the map $g(\cdot)$ is also order-preserving on $\mathbb{S}^n_{+}$ with respect to the Loewner partial order.
\end{itemize}

\smallskip
\paragraph{Comparison of various gradient-based surrogates}
Applying the Gaussian Poincar\'e inequality (see \cite{bakry2013analysis} and references therein), the covariance $\Cov_\refparams(\mathcal{J}(\design_k, \refparams))$ satisfies\footnote{This is a direct consequence of the Poincar\'e inequality for univariate functions by considering, for any direction $\bvec{v} \in \mathbb{R}^{\Nd}$ with $0 < \|\bvec{v}\| < \infty$, we have
\(
\bvec{v}^\top \Cov(\mathcal{J}) \bvec{v} = \mathrm{var}( \mathcal{J}^\top \bvec{v} )\leq \Expect{}{ \| \nabla(\mathcal{J}^\top \bvec{v})\|^2 } = \bvec{v}^\top \Expect{}{\nabla \mathcal{J}  \nabla \mathcal{J}^\top} \bvec{v}.
\)
}
\begin{equation}
\Cov_\refparams\big(\mathcal{J}(\design_k, \refparams)\big)\preceq \int  \Jac_{\refparams} \mathcal{J}(\design_k, \refparams) \Jac_\refparams \mathcal{J}(\design_k, \refparams)^{\top} \refD(\refparams) \mathrm{d} \refparams =: \Hmat_{\mathrm{G}}(\design_k, \mathcal{T}^{k{-}1}) \in \mathbb{R}^{\Nd \times \Nd}.
\label{eq:cov_bound}
\end{equation}
Using the first identity above, the covariance-based surrogate \eqref{eq:cov_bound} is bounded by
\begin{equation}
\EIG^{\mathrm{C}}_k(\design_k) \leq \frac{1}{2}\log \det \big(  \I{\Nd} + \Cnoise(\design_k)^{-1}\Hmat_{\mathrm{G}}(\design_k, \mathcal{T}^{k{-}1}) \big) \eqqcolon \EIG^{\mathrm{G}}_k(\design_k).
\label{eq:cov_grad_bound}
\end{equation}
This leads to a comparison of our bounds with the work of \cite[Section 4.1]{ChenArnaudBaptistaZahm:2025:1}, which uses the bounding matrix $\Hmat_{\mathrm{G}}(\design_k, \mathcal{T}^{k{-}1})$ to derive an EIG surrogate for batched sensor selection with i.i.d. Gaussian observation noises.
Translating to our sequential setting, the setup of \cite[Section 4.1]{ChenArnaudBaptistaZahm:2025:1} gives the noise covariance $\Cnoise(\design_k) = \sigma^2 \I{\Nd}$, the bounding matrix
\[
\Hmat_{\mathrm{G}}(\design_k, \mathcal{T}^{k{-}1}) = \bvec{W}(\design_k) \left(\int \Jac_\refparams ( \mathcal{F}\circ \cT^{k{-}1}(\refparams) ) \Jac_\refparams ( \mathcal{F}\circ \cT^{k{-}1}(\refparams) )^\top \refD(\refparams) \mathrm{d}\refparams\right) \bvec{W}(\design_k),
\]
and the corresponding design criterion
\begin{equation}
\design_k^* \in \Argmax_{\design_k} \EIG_k^{\prime}(\design_k) := \Argmax_{\design_k} \frac12 \mathrm{trace} \left( \Cnoise(\design_k)^{-1} \Hmat_{\mathrm{G}}(\design_k,\cT^{k{-}1}) \right)  .
\label{eq:chen_etal}
\end{equation}
Applying the second identity above, the iEIG, the covariance-based surrogate, its gradient-based bound, and the trace surrogate satisfy the ordering
\[
\EIG_k(\design_k) \leq \EIG_k^{\mathrm{C}}(\design_k)  \leq \EIG_k^{\mathrm{G}}(\design_k)  < \EIG_k^{\prime}(\design_k),
\]
The gap between $\EIG_k^{\mathrm{G}}(\design_k)$ and $\EIG_k^{\prime}(\design_k)$ increases as $\Hmat_{\mathrm{G}}$ grows, i.e., as more information is acquired from the data.
This also allows a comparison between the trace surrogate $\EIG_k^{\prime}(\design_k)$ and our information-based surrogate $\EIG_k^{\mathrm{I}}(\design_k)$. By linearity of the trace and expectation, we have
\(
 \mathrm{trace}\big(\mathcal{H}_{\mathrm{I}}(\design_k)\big) = \mathrm{trace} \big( \Cnoise(\design_k)^{-1} \Hmat_{\mathrm{G}}(\design_k,\cT^{k{-}1}) \big)
\)
for additive Gaussian noise. Consequently, the iEIG, the information-based surrogate, and the trace surrogate satisfy the ordering
\[
\EIG_k(\design_k) \leq \EIG_k^{\mathrm{I}}(\design_k) < \EIG_k^{\prime}(\design_k).
\]
Similarly, the gap between $\EIG_k^{\mathrm{I}}(\design_k)$ and $\EIG_k^{\prime}(\design_k)$ increases as $\Hmat_{\mathrm{I}}$ grows.

\smallskip
\paragraph{Comparison of information- and covariance-based bounds}
Although the alternative gradient-based bound $\EIG^{\mathrm{G}}_k(\design_k)$ in \eqref{eq:cov_grad_bound} and the information-based bound $\EIG^{\mathrm{I}}_k(\design_k)$ (cf.~\eqref{eq:iEIG_uB} and \eqref{eq:Fisher_gauss}) involve the same ingredients, neither a general ordering nor a sharp comparison is currently available in multidimensional, nonlinear settings. Moreover, the relationship between the covariance-based bound and the information-based bound is also unclear.
We explore this relationship numerically in the Supplementary Material using a toy problem from~\cite[Section 5.1]{HuanMarzouk:2013:1}, and in a more practical example in~\cref{sec:num1}. These results suggest that the covariance-based bound tends to be tighter than the information-based bound, at least when $\Nd<\Nm$.
However, the covariance-based bound relies on the assumption of additive Gaussian noise. While it is currently unclear how to extend this bound to more general likelihood functions, the gradient-based bound depends only on the Fisher information, which is well-defined for a broad class of likelihoods. Therefore, in the following sections, we focus on the gradient-based bound and explore its use for scalable sequential design.

\smallskip
\paragraph{On intrusive design settings}
For intrusive designs where the design may change the governing equations, e.g., finding the optimal functional form of source terms, fixing a candidate set is generally not appropriate. In such cases, an optimization-based approach is needed, which typically requires re-evaluating the forward model at new designs for each sample. This contrasts with the fixed candidate setting, where the upper bound can be evaluated at \emph{all} candidate designs using just $N$ forward solves (and at most $N\Nd$ adjoint solves when using the gradient-based bound) in a single Monte Carlo approximation. However, the bounds proposed here can still serve as fast surrogates in intrusive settings.

\subsection{Scalable amortized inference}\label{subsec:sOED_composite}
The goal of amortized inference is to construct a conditional map $\cS^k_{\params \given \data} \colon\mathbb{R}^{\Nd} \times \mathbb{R}^{\Nm} \rightarrow \mathbb{R}^{\Nm}$, which approximately couples the stage-$k$ posterior to a reference density for any observed data. This map can be built during experimentation, before data are collected, and then enables rapid online inference of the posterior $\target(\params\mid\Hist_k)$ for newly observed data and derive the iEIG bound at the next stage. We combine the data-free LIS with tensor trains to achieve this.

\smallskip
\paragraph{Conditional map}
At stage $k$, suppose we have a transport map $\cT^{k-1}$---either the previous conditional map or other forms of transport maps---and a standard multivariate Gaussian reference $\refD$ such that $\widehat{\target}(\params\given \Hist_{k{-}1}):=\pushforward{{(\cT^{k{-}1})}} \refD(\params)$ approximates the stage-$k$ prior $\target(\params \given \Hist_{k{-}1})$. After selecting the optimal design $\design_k^*$, we construct the conditional map as follows.

\begin{enumerate}[leftmargin=19pt]
    \item Define a joint reference density $\refD(\refparams_y,\refparams_m) = \refD(\refparams_y)\refD(\refparams_m)$ and a \emph{joint precondition map},
    $\cS^{k{-}1}_{\data,\params}= [\refparams_y,{\cT^{k{-}1}}(\refparams_m)]^\top$,
    such that $\pushforward{(\cS^{k{-}1})}\refD(\data,\params) = \refD(\data)\,\widehat{\target}(\params \given \Hist_{k{-}1})$.
    \item Pull back the stage-$k$ joint density for the data and parameter random variables,
\begin{equation*}
\target(\data_k,\params\given \design_k^*,\Hist_{k{-}1}) = \likelihood{\data_k}{\params,\design_k^*}\,\target(\params\given\Hist_{k{-}1}),
\end{equation*}
under the preconditioning map, $\cS^{k{-}1}_{\data,\params}$, which yields
\begin{align}
\pullback{({\cS^{k{-}1}_{\data,\params}})}\target(\ratioparams_y,\ratioparams_{m}\given \design_k^*,\Hist_{k{-}1}) &= \frac{\likelihood{\ratioparams_y}{\cT^{k{-}1}(\ratioparams_m),\design_k^*}\,\target(\cT^{k{-}1}(\ratioparams_m)\given \Hist_{k{-}1})}{\refD(\ratioparams_y)\,\widehat{\target}(\cT^{k{-}1}(\ratioparams_m)\given \Hist_{k{-}1})}\refD(\ratioparams_y,\ratioparams_m) \nonumber \\
&\approx \likelihood{\ratioparams_y}{\cT^{k{-}1}(\ratioparams_m),\design_k^*}\,\refD(\ratioparams_m) \eqqcolon q^k_{\data,\params}(\ratioparams_y,\ratioparams_m). \label{eq:qMap}
\end{align}
\item Approximate the joint density $q^k_{\data,\params}(\ratioparams_y,\ratioparams_m)$ using the deep approximation as in~\cref{subsec:DIRT} to obtain a lower-triangular \emph{incremental joint map} 
\begin{equation}
\mathcal{Q}_{\data,\params}^k (\refparams_y,\refparams_m) = \begin{bmatrix*}[l]
		\mathcal{Q}^k_{\data}(\refparams_y)\vspace{2pt}
		\\
		\mathcal{Q}^k_{\params \given \data}(\refparams_y, \refparams_m)
	\end{bmatrix*}
    =
    \begin{bmatrix*}[l]
		\ratioparams_y
		\\
	   \ratioparams_m
	\end{bmatrix*},\label{eq:incrementalMap}
\end{equation}
that approximately pushes forward the reference density $\refD(\refparams_y,\refparams_m)$ to $q^k_{\data,\params}(\ratioparams_y,\ratioparams_m)$.
\item The composite map, $\cS^k_{\data,\params} = \cS^{k{-}1}_{\data,\params} \circ \cQ^{k}_{\data,\params} $, which also takes a lower triangular form of
\[
\cS_{\data,\params}^k (\refparams_y,\refparams_m) = \begin{bmatrix*}[l]
		\;\cQ^k_{\data}(\refparams_y)\vspace{2pt}
		\\
		(\cT^{k{-}1} \circ \cQ^k_{\params \given \data})(\refparams_y, \refparams_m)
	\end{bmatrix*}
    =
    \begin{bmatrix*}[l]
		\data_k
		\\
	    \params
	\end{bmatrix*},\label{eq:jointMap}
\]
approximates the joint density $\target(\data_k,\params\given \design_k^*,\Hist_{k{-}1})$. The resulting \emph{conditional map}
\begin{equation}
\cS^{k}_{\params \given \data}(\data_k, \refparams_m) = \left(\cT^{k{-}1} \circ \cQ^k_{\params \given \data}\right) \left((\cQ^k_{\data})^{-1}(\data_k), \refparams_m\right),
    \label{eq:conditionalMap}
\end{equation}
approximately pushes forward the reference density $\refD(\refparams_m)$ to the  conditional marginal density $\target(\params\given \data_k,\design_k^*,\Hist_{k{-}1})$, which is the posterior, for any data $\data_k$ observed at stage $k$.
\end{enumerate}

\smallskip
Once an observed data $\data_k$ is collected at $k$-th experimental stage, we define $\cT^{k}(\refparams_m) := \cS^{k}_{\params \given \data}(\data_k, \refparams_m)$ to approximate the updated posterior. It also serves as an input to~\cref{Alg:DLIS_combined} for finding optimal experimental conditions in the next stage.

\smallskip
\paragraph{Subspace accelerated conditional map}

When building the incremental joint map $\cQ^{k}_{\data,\params}$ described above, the number of density evaluations scales linearly with the joint data--parameter dimension and quadratically with the tensor ranks, which may grow with dimension. To overcome this scalability bottleneck, we embed the conditional map construction within the LIS that naturally comes with our information-based iEIG surrogate.

For the transformed parameters $(\ratioparams_y,\ratioparams_m)$, the joint density $q^k_{\data,\params}$ defined in \eqref{eq:qMap} is equipped with a Gaussian parameter prior $\refD(\ratioparams_m)$ and the likelihood $\likelihood{\ratioparams_y}{\cT^{k{-}1}(\ratioparams_m),\design_k^*}$. Let $\LISDF_{k}$ be the $\LISdim_k$-dimensional data-free LIS basis that captures directions along which the transformed parameters $\ratioparams_m$ are most informed by the likelihood. Decomposing the parameter $\ratioparams_m = \ratioparams^{r_k}_{m}+\ratioparams_m^{\perp_k}$ into its likelihood-informed components ($\ratioparams^{r_k}_{m})$ and uninformed components ($\ratioparams_m^{\perp_k}$) and following~\eqref{eq:li_posterior}, we can approximate the joint density $q^k_{\data,\params}(\ratioparams_y,\ratioparams_m)$ as
\begin{equation}
\begin{split}
q^k_{\data,\params}(\ratioparams_y,\ratioparams_m) \approx \widetilde{q}^k_{\data,\params}(\ratioparams_y,\ratioparams_m;\Proj_{k})
&= \underbrace{\Rlikelihood{\ratioparams_y}{\mathcal{T}^{k{-}1}(\ratioparams_m^{r_k}),\design_k^*}\,\refD(\ratioparams_m^{r_k})}_{\widetilde{q}^k_{\data,\params}(\ratioparams_y,\ratioparams_m^{r_k})}\,\refD(\ratioparams_m^{\perp_k}\given \ratioparams_m^{r_k}).
\end{split}
\label{eq:approx_intermediat}
\end{equation}
We obtain the data-free LIS basis $\LISDF_k$ from the leading eigenvectors of the average Fisher information matrix $\overline{\Hmat}_\mathrm{I}(\design_k^*,{\cT^{k{-}1}})$ computed in~\cref{Alg:DLIS_combined}. Using the bound on the expected KL divergence~\eqref{eq:EIG_bound}, we can choose the LIS dimension $\LISdim_k$ as the smallest integer $\LISdim$ such that
\begin{equation}
\frac{1}{2}\sum_{j=\LISdim+1}^{\Nm} \log\left( 1 + \lambda_j\left(\overline{\Hmat}_\mathrm{I}(\design_k^*,\mathcal{T}^{k{-}1})\right)\right) \leq \epsilon_I,
\end{equation}
for some error tolerance $\epsilon_I$ on the approximation to the joint density $q^k_{\data,\params}$.
Once the data-free LIS is determined, let $\widetilde{\ratioparams}_m^k = \LISDF_k^\top\ratioparams_m \in \mathbb{R}^{\LISdim_k}$ be the coefficient associated with the LIS basis. Then one can approximate the reduced-dimensional joint density, $\widetilde{q}^k_{\data,\params}(\ratioparams_y,\ratioparams_m^{r_k}) \equiv \widetilde{q}^k_{\data,\params}(\ratioparams_y,\widetilde{\ratioparams}_m^k)$, to build a smaller KR map $\widetilde{\cQ}^k_{\data,\params}:\mathbb{R}^{\Nd+\LISdim_{k}} \rightarrow \mathbb{R}^{\Nd+\LISdim_k}$,
\[
\widetilde{\cQ}^k_{\data,\params}(\refparams_y,\refparams_m^k) = [\widetilde{\cQ}^k_{\data}(\refparams_y), \widetilde{\cQ}^k_{{\params}\given \data}(\refparams_y,\refparams_m^k)]^\top, \quad \refparams_m^k \in \mathbb{R}^{\LISdim_k}
\]
such that
\(
\pushforward{({\widetilde{\cQ}^k_{\data,\params}})}\refD(\ratioparams_y,\widetilde{\ratioparams}_{m}^k) \approx \widetilde{q}^k_{\data,\params}(\ratioparams_y,\widetilde{\ratioparams}_{m}^k).
\)
Embedding the map $\widetilde{\cQ}^k_{\data,\params}$ into a linear map defined using $\LISDF_k$, we obtain the incremental map on the full space
\begin{align}
    \cQ^{k}_{\data,\params}(\refparams_y,\refparams_m) &= \begin{bmatrix*}[l]
		\widetilde{\cQ}^{k}_{\data}(\refparams_y) \vspace{2pt}
		\\
		{\cQ}^{k}_{\params \given \data}(\refparams_y,{\refparams_m})
	\end{bmatrix*},
    \label{eq:ratio_fullSpace}
\end{align}
where
\(
    {\cQ}^{k}_{\params \given \data}(\refparams_y,{\refparams_m}) = \LISDF_{k}\widetilde{\cQ}^{k}_{\params \given \data}(\refparams_y,\LISDF_{k}^\top\refparams_m)+(\I{\Nm}-\LISDF_{k}\LISDF_{k}^\top)\refparams_m.
\)

\subsection{Overall algorithm with restart}\label{subsec:sOED_restart}

The method described so far can be effective for sOED problems with a small number of experimental stages.
However, several practical challenges arise when dealing with a large number of experimental stages. First, the incremental accumulation of error, due to the recursive approximation, can make it computationally expensive to maintain a sufficiently small approximation error across the experimental stages.
This challenge is compounded by the use of conditional transport maps, as they can only probabilistically guarantee a sufficiently accurate approximation to the posterior~\cite[Appendix A.2]{CuiDolgovZahm:2023:1}.
Second, the continuous addition of layers to the composite map with experimental stages increases the computational cost of sampling.
Third, relying solely on the data-free LIS may lead to an overestimation of the number of important directions.

To address these issues, we incorporate a ``restart'' step into the conditional map construction. Specifically, every $\ell$ steps---after collecting data from previous experiments but before designing the experimental conditions for the next stage---we discard the previous conditional map $\cS^\ell_{\params\given\data}$ and rebuild a map $\cT^\ell$ from scratch by targeting the stage-$\ell$ posterior directly. This process produces a more accurate approximation of the current posterior. To ensure scalability, the approximation begins by constructing the data-dependent LIS, followed by the reduced-dimensional posterior approximation. Specifically, we define an importance sampling estimator using the previous conditional map $\cS^\ell_{\params\given\data}$ to compute the Gram matrix~\eqref{eq:gram} for building a data-dependent LIS. This procedure is summarized in~\Cref{Alg:DTLIS}.

\begin{algorithm}[h]
\caption{Construction of the data-dependent LIS basis, $\LISDD$, for a posterior $\target(\params \given \data) \propto \targetUnnorm(\params\given\data) = \likelihood{\data}{\params}\refD(\params)$. Inputs are a transport map $\cS_{\params \given \data}$ such that $\pushforward{({\cS_{\params \given \data}})}\refD(\params) = \hat{\target}(\params \given \data) \approx \target(\params \given \data)$, a truncation tolerance $\epsilon_{D}$, and a sample sizes $N$. \label{Alg:DTLIS}}
\begin{algorithmic}[1]
\Function{DDLIS}{$\cS_{\params \given \data},\epsilon_{D},N$}
        \State Draw sample $\refparams^{(i)}\sim \refD$ and evaluate $\params^{(i)} = {\cS_{\params \given \data}}(\refparams^{(i)})$ for $i = 1, \ldots, N$
        \State Compute the Gram matrix using self-normalized importance sampling \vspace{-12pt}$$\overline{\Hmat}(\bvec{y}) = \textstyle \frac{1}{\sum_{i=1} W_i}\sum_{i=1}^N W_{i}\nabla_\params\log\likelihood{\data}{\params^{(i)}}\nabla_\params\log\likelihood{\data}{\params^{(i)}}^{\top}, \quad W_{i} = \frac{\targetUnnorm(\params^{(i)}|\data)}{\hat{\target}(\params^{(i)}|\data)} \vspace{-12pt}$$
        \State Compute eigendecomposition $\overline{\Hmat}(\bvec{y}) = \LISDD \bvecS{\Lambda} \LISDD^\top$
        \State Set $\LISDD = \LISDD(:,1:\LISdim)$, with $\LISdim$ chosen using~\eqref{eq:DI_bound} to ensure $\frac{\sqrt{\kappa}}{2}\sqrt{\sum_{k=r+1}^{\Nm}\lambda_k(\overline{\Hmat}(\bvec{y})))} \leq \epsilon_{D}$
        \State \Return $\LISDD$
\EndFunction
\end{algorithmic}
\end{algorithm}

\begin{figure}
\centering
\begin{tikzpicture}
\node[inner sep=0pt] (a) at (0,0)
{\includegraphics[width=\textwidth]{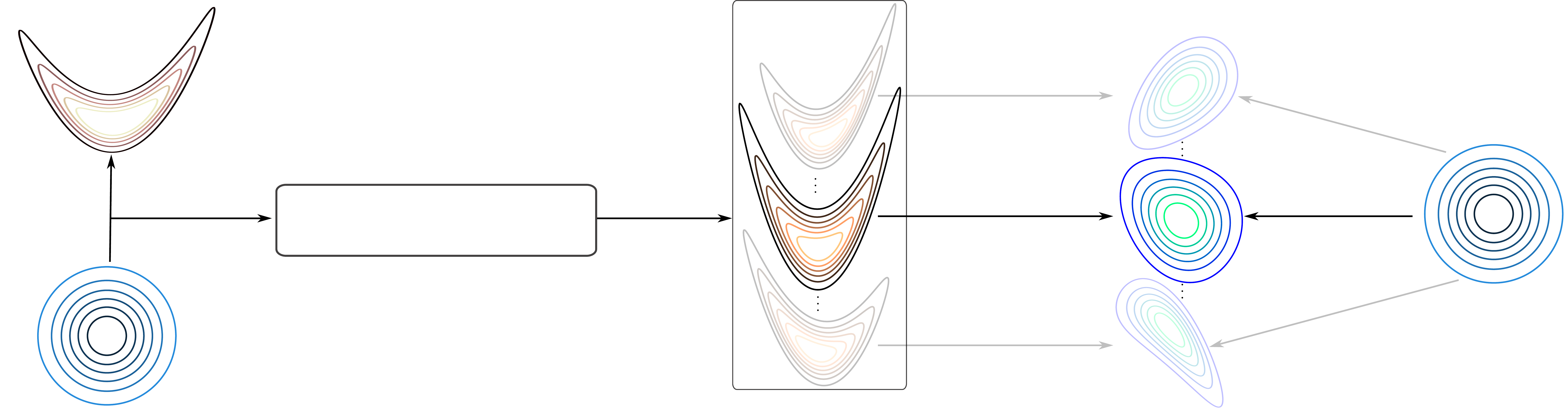}};
\node at (-6.65,1.8) {\small{$\pi_{\Hist_{k-1}}$}};
\node at (-6.7,-2.25) {\small{$\refD$}};
\node at (0.3,2.35) {\small{$\pi_{\data_k,\params \given \design^{*}_k,\Hist_{k-1}}$}};
\node at (0.3,-2.1) {\small{amortized inference}};
\node at (-3.4,0.4) {\small{iEIG surrogate}};
\node at (-3.4,-0.15) {\small{$\argmax \EIG^{\mathrm{X}}_k(\design_k)$}};
\node at (-1.3,0.2) {\small{$\design^{*}_k$}};
\node at (-6.1,0.2) {\small{$\cT^{k-1}$}};
\node at (2.2,0.3) {\small{$\cT^{k-1}$}};
\node at (5.5,0.3) {\small{$\cQ^{k}$}};
\node at (7.05,0.8) {\small{$\refD$}};
\label{fig:alg_workflow}
\end{tikzpicture}\vspace{-25pt}
\caption{The workflow of one experimental stage ($k$) of~\cref{Alg:sOED_wRestart}.}
\end{figure}

In each restart, the procedure for constructing the transport map with the data-dependent LIS closely follows the subspace acceleration approach described in the previous subsection.
The key differences are that the LIS dimension is determined by the data-dependent bound~\eqref{eq:DI_bound} and a non-conditional transport map is built by directly approximating the fixed posterior.
Assuming the stage $\ell$ requires a restart---\ie, an accurate approximation to the stage $\ell{-}1$ posterior is needed---the transport map using the restart is outlined as follows.
Let $\LISDD_{\ell{-}1} \in \mathbb{R}^{\Nm \times \LISdim_{\ell{-}1}}$ denote a basis of the data-dependent LIS for the current posterior $\target(\params \given \Hist_{\ell{-}1})$, computed using~\cref{Alg:DTLIS}. The parameter $\params$ can be decomposed as $\params = \params^{\LISdim_{\ell{-}1}} + \params^{\perp_{\ell{-}1}}$. Let $\widetilde{\params}^{\ell{-}1} \in \mathbb{R}^{\LISdim_{\ell{-}1}}$ be the coefficients of the reduced parameter $\params^{\LISdim_{\ell{-}1}}$ associated with $\LISDD_{\ell{-}1}$. One can first construct a lower-dimensional map $\widetilde{\cT}^{\ell{-}1}\colon\mathbb{R}^{\LISdim_{\ell{-}1}} \rightarrow \mathbb{R}^{\LISdim_{\ell{-}1}}$ satisfying $\pushforward{(\widetilde{\cT}^{\ell{-}1})}\refD(\widetilde{\params}^{\ell{-}1}) \approx \widetilde{\target}(\LISDD_{\ell{-}1}\widetilde{\params}^{\ell{-}1} \given \Hist_{\ell{-}1})$, with the reduced-dimensional posterior $\widetilde{\target}(\cdot \given \Hist_{\ell{-}1})$ defined in~\eqref{eq:li_posterior}. Then, embedding the smaller map $\widetilde{\cT}^{\ell{-}1}$ into a linear map as in~\cref{subsec:sOED_composite}, we obtain the full-dimensional transport map
\begin{equation}
\cT^{\ell{-}1}(\refparams) = \LISDD_{\ell{-}1}^{}\widetilde{\cT}^{\ell{-}1}(\LISDD_{\ell{-}1}^\top\refparams)+\big(\I{\Nm}-\LISDD_{\ell{-}1}^{}\LISDD_{\ell{-}1}^\top\big)\refparams, \quad \refparams \in \mathbb{R}^{\Nm},\label{eq:restartMap}
\end{equation}
which approximately couples $\refD = \mathcal{N}(\bvec{0},\I{\Nm})$ and ${\target}(\params \given \Hist_{\ell{-}1})$. The new map $\cT^{\ell{-}1}$ can then be used in~\cref{Alg:DLIS_combined} for finding new optimal designs and building the data-free LIS for the amortized inference in the new experimental stage. The overall algorithm is presented in~\cref{Alg:sOED_wRestart}, and the workflow for one stage is visualized in~\cref{fig:alg_workflow}.

\begin{algorithm}[h!]
	\caption{Greedy sOED with restart. Inputs are the number of experimental stages $K$, truncation tolerances $\epsilon_D$ and $\epsilon_I$, a stacked forward map $\mathcal{F}$, and sample sizes $N_r$, $N$. Without loss of generality, we set the initial prior $\refD \sim \mathcal{N}(\bvec{0},\I{\Nm})$.}
	\label{Alg:sOED_wRestart}
	\begin{algorithmic}[1]
		\State Initialize $\cT^0 = \mathcal{I}$
		\For{k=1, \ldots, K}
        \If{k $>$ 1}
        \State Get new data $\data_{k{-}1}$ and define $\cT^{k{-}1}(\cdot) = \cS^{k{-}1}_{\params \given \data}(\data_{k{-}1}, \cdot)$ for posterior estimation
        \If{restart}
        \State $\LISDD_{k{-}1} \leftarrow \mathrm{DDLIS}(\cS^{k{-}1}_{\params \given \data},\epsilon_D,N_r)$ using~\cref{Alg:DTLIS}
        \State Build the reduced map $\widetilde{\cT}_{k{-}1}$ such that $\pushforward{({\widetilde{\cT}_{k{-}1}})}\refD(\widetilde{\params}^{k{-}1}) \approx \widetilde{\target}(\LISDD_{k{-}1}\widetilde{\params}^{k{-}1} \given \Hist_{k{-}1})$
        \State Build the full-dimensional map ${\cT^{k{-}1}}(\cdot)$ as in~\eqref{eq:restartMap} for posterior estimation
        \EndIf
        \EndIf
        \State $\design_k^*,\LISDF_{k} \leftarrow \mathrm{iEIGUB}({\cT^{k{-}1}},\mathcal{F},\epsilon_I,N)$ using~\cref{Alg:DLIS_combined}
        \State  Define the joint preconditioning map $\cS^{k{-}1}_{\data,\params}= [\refparams_y,{\cT^{k{-}1}}(\refparams_m)]^\top$
        \State Build the incremental joint map $\cQ^k_{\data,\params}$ as in~\eqref{eq:incrementalMap}
        \State Define $\cS_{\data,\params}^{k} = \cS^{k{-}1}_{\data,\params} \circ \cQ^{k}_{\data,\params}$ and extract the new conditional map $\cS^{k}_{\params \given \data}$ as in~\eqref{eq:conditionalMap}
		\EndFor
	\end{algorithmic}
\end{algorithm}

The additional computational effort of the restart step is justified by the improved accuracy of posterior approximations. This, in turn, can enhance the stability of the approximations and potentially lower the tensor ranks in subsequent experimental stages. Furthermore, the computational cost of constructing transport maps using tensor train decompositions can be significantly reduced by utilizing surrogate models for the forward map. These surrogates can be built using neural networks or, as demonstrated in our numerical experiments, through proper orthogonal decomposition and the discrete empirical interpolation method~\cite{ChaturantabutSorensen:2010:1}.

We summarize the computational cost of the main steps of~\cref{Alg:sOED_wRestart}. At each stage, the dominant cost arises from forward model evaluations, occurring in the evaluation of the iEIG upper bound (line 9), as well as the construction of transport maps (the conditional map in line~11, and the non-conditional maps in lines~6--7 when using restarts). We outline below where these expensive solves are replaced by a cheaper surrogate in our numerical experiments.
\begin{enumerate}[label=(\arabic*), leftmargin=19pt]
\item \textbf{Evaluation of iEIG upper bound:} Requires $N$ forward model evaluations. If a gradient-based bound is used, at most $N \Ne \Nd$ adjoint solves are also required.
\item \textbf{Conditional map construction:} Let $M_k$ and $R_k$ denote the maximum basis size and rank over the $L_k$ layers of the deep tensor construction described in~\cref{subsec:DIRT}. This step requires $O((\Nd+r_k)L_kM_kR_k^2)$ evaluations of the unnormalized joint density $\widetilde{q}^k_{\data,\params}$. In our simulations, the forward map in the likelihood of~\eqref{eq:approx_intermediat} is replaced by a surrogate.
\item \textbf{Non-conditional map construction:} Constructing the data-dependent LIS requires $N_r$ forward model evaluations together with at most $N_r\Nd$ adjoint solves. Using the same notation $r_k$, $R_k$, and $M_k$ as in~(2), the subsequent transport map construction requires $O(r_k L_k M_k R_k^2)$ evaluations of the unnormalized stage-$k$ posterior density, where the forward map is replaced by a surrogate in our simulations.
\end{enumerate}

\section{Numerical results} We demonstrate the proposed approach on two model problems.

\subsection{Problem 1: optimal sensor selection for diffusivity field estimation}\label{sec:num1}
We first consider an inverse problem governed by an elliptic PDE commonly used to model the flow of a fluid through a porous medium. The PDE is defined in $\Omega = [0,1]^2$ with left, right, top and bottom boundaries denoted as $\Gamma_L, \Gamma_R, \Gamma_T, \Gamma_B$, respectively:
\begin{align*}
-\nabla \cdot \left(\kappa(m) \nabla u \right) &= 0 && \hspace{-6em} \text{ in } \Omega, \\
\kappa(m)\nabla u \cdot n &= 0  && \hspace{-6em}  \text{ on } \Gamma_T \cup \Gamma_B, \\
u &= 1+y/2 && \hspace{-6em}  \text{ on } \Gamma_L, \\
u &= -\sin(2\target y)-1 && \hspace{-6em}  \text{ on } \Gamma_R.
\end{align*}
We consider the inverse problem of inferring the log-diffusivity field $m = \log(\kappa)$ from pressure $u$ measured at finitely many locations in $\Omega$.
We assume $m$ follows a Gaussian prior  $m\sim\mathcal{N}(0,C_{\mathrm{pr}})$ with $C_{\text{pr}}$ defined by the covariance kernel $\exp(-\frac{1}{2\ell^2} \norm{}{x-z}^{2}),$ $\ell = \frac{1}{\sqrt{50}}$.

\begin{figure}
\vspace{-24pt}
\centering
\begin{tikzpicture}
\node[inner sep=0pt] (a) at (0,0)
{\includegraphics[width=0.7\textwidth, trim={0 5em 0 5em}, clip]{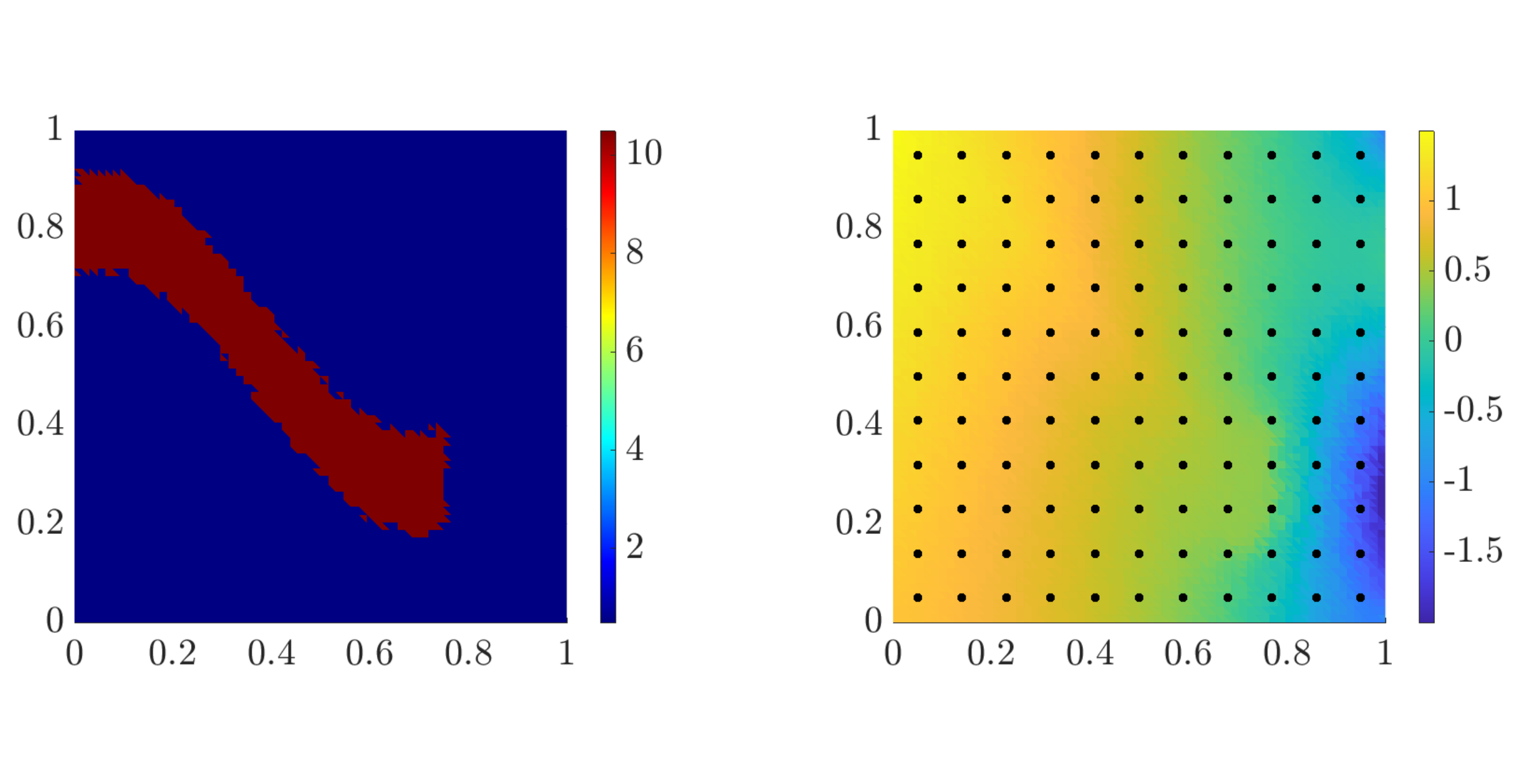}};
\node at (-3.3,2.1) {{$\kappa_{\text{true}}$}};
\node at (2.8,2.1) {{$u(\kappa_{\text{true}})$}};
\label{fig:Darcy_setup}
\end{tikzpicture}\vspace{-15pt}
\caption{The ``true'' diffusivity field used to synthesize data (left) and the corresponding pressure field $u$ (right). The $\Nd = 121$ candidate locations for the design problem are visualized as black dots in the right figure.}
\end{figure}

For the sequential design problem, we fix $\Nd = 121$ equally-spaced candidate locations for measuring the pressure $u$ and assume the measurement at each candidate location is corrupted by independent mean-zero Gaussian noise with standard deviation $0.2$. In each experiment, we choose one location from this set that maximizes the incremental EIG bound~\eqref{eq:EIG_k}. The true diffusivity field used to synthesize the data at each experimental stage, as well as the corresponding PDE solution $u$ and the candidate sensor locations are visualized in~\cref{fig:Darcy_setup}.

After finite element discretization with second-order Lagrange elements on a triangular mesh of size $h = \frac{1}{32}$ in each coordinate direction, the unknown log-diffusivity field is characterized by a vector of nodal coefficients $\bvec{m} \in \mathbb{R}^{\Nm}$ with $\Nm = 4225$. The forward solves as well as the adjoint-based construction of the Jacobians needed for the upper bound on the iEIG~\eqref{eq:iEIG_pullback} are performed in \texttt{FastFins}~\cite{Cui:2022}.

\smallskip
\paragraph{Optimal designs and comparisons}
We use eight design stages, each selecting an optimal location for measuring $u$ using algorithm~\cref{Alg:sOED_wRestart} with $N=100$ samples and truncation tolerances of $\epsilon_G = 0.01, \epsilon_I = 0.02$, and transport maps constructed using the \texttt{deep-transport} package~\cite{Cui:2023}.
To speed up the computation of the tensor train density approximations, we build a reduced-order model (ROM) restricted to the data-free LIS at all the candidate sensor locations.
Specifically, we use the discrete empirical interpolation approach combined with proper orthogonal decomposition as in~\cite{ChaturantabutSorensen:2010:1}, with $1000$ samples drawn from the LIS to create the snapshots.
The LIS at all the candidate designs that satisfies~\eqref{eq:EIG_bound} with a more conservative tolerance of $0.01$ has dimension $71$, and solving the resulting ROM, which has a relative $L_2$ error of approximately $0.003$, is approximately 40 times faster than solving the full model.
The ROM is only used to speed up the construction of the transport maps in each design stage. To avoid ROM error perturbing the design selection, the average information matrix appearing in the iEIG bound~\eqref{eq:iEIG_uB} is computed using the high-fidelity model.
The computational cost of choosing the optimal designs and approximating the posteriors, along with the approximation errors, is presented in~\cref{table:cost_darcy} for both the covariance- and information-based bounds.

\begin{table}
\footnotesize
\caption{The dimensions of the data-free and data-dependent LISs for each stage (top and bottom rows of column 2, respectively), as well as the Hellinger distance between the true and approximate stage-$k$ posterior (column 3) and the effective sample size per sample ($\text{ESS}/N$, column 4). The last column lists the total number of ROM solves required to approximate the posterior in each stage.} 
\centering
\begin{tabular}{|c|c|c|c|c|}
\hline
\textbf{$k$}              & \textbf{\# LIS} & \textbf{$\mathcal{D}_{\text{H}}$} & \textbf{ESS/N}         & \textbf{\# ROMs}         \\ \hline
\multirow{2}{*}{\textbf{1}} & 15              & \multirow{2}{*}{0.127}            & \multirow{2}{*}{0.788} & \multirow{2}{*}{108,500} \\ \cline{2-2}
                            & 28              &                                   &                        &                          \\ \hline
\multirow{2}{*}{\textbf{2}} & 14              & \multirow{2}{*}{0.211}            & \multirow{2}{*}{0.523} & \multirow{2}{*}{111,631} \\ \cline{2-2}
                            & 31              &                                   &                        &                          \\ \hline
\multirow{2}{*}{\textbf{3}} & 13              & \multirow{2}{*}{0.176}            & \multirow{2}{*}{0.752} & \multirow{2}{*}{210,552} \\ \cline{2-2}
                            & 35              &                                   &                        &                          \\ \hline
\multirow{2}{*}{\textbf{4}} & 16              & \multirow{2}{*}{0.253}            & \multirow{2}{*}{0.362} & \multirow{2}{*}{244,807} \\ \cline{2-2}
                            & 40              &                                   &                        &                          \\ \hline
\multirow{2}{*}{\textbf{5}} & 12              & \multirow{2}{*}{0.274}            & \multirow{2}{*}{0.133} & \multirow{2}{*}{426,157} \\ \cline{2-2}
                            & 42              &                                   &                        &                          \\ \hline
\multirow{2}{*}{\textbf{6}} & 13              & \multirow{2}{*}{0.231}            & \multirow{2}{*}{0.662} & \multirow{2}{*}{443,021} \\ \cline{2-2}
                            & 43              &                                   &                        &                          \\ \hline
\multirow{2}{*}{\textbf{7}} & 15              & \multirow{2}{*}{0.249}            & \multirow{2}{*}{0.606} & \multirow{2}{*}{744,992} \\ \cline{2-2}
                            & 44              &                                   &                        &                          \\ \hline
\multirow{2}{*}{\textbf{8}} & 15              & \multirow{2}{*}{0.266}            & \multirow{2}{*}{0.542} & \multirow{2}{*}{762,755} \\ \cline{2-2}
                            & 45              &                                   &                        &                          \\ \hline
\end{tabular}\label{table:cost_darcy}
\quad
\begin{tabular}{|c|c|c|c|c|}
\hline
\textbf{$k$}              & \textbf{\# LIS} & \textbf{$\mathcal{D}_{\text{H}}$} & \textbf{ESS/N}         & \textbf{\# ROMs}         \\ \hline
\multirow{2}{*}{\textbf{1}} & 16              & \multirow{2}{*}{0.110}            & \multirow{2}{*}{0.898} & \multirow{2}{*}{121,365} \\ \cline{2-2}
                            & 32              &                                   &                        &                          \\ \hline
\multirow{2}{*}{\textbf{2}} & 11              & \multirow{2}{*}{0.151}            & \multirow{2}{*}{0.802} & \multirow{2}{*}{102,889} \\ \cline{2-2}
                            & 31              &                                   &                        &                          \\ \hline
\multirow{2}{*}{\textbf{3}} & 14              & \multirow{2}{*}{0.180}            & \multirow{2}{*}{0.755} & \multirow{2}{*}{225,091} \\ \cline{2-2}
                            & 37              &                                   &                        &                           \\ \hline
\multirow{2}{*}{\textbf{4}} & 12              & \multirow{2}{*}{0.200}            & \multirow{2}{*}{0.728} & \multirow{2}{*}{218,457} \\ \cline{2-2}
                            & 37              &                                   &                        &                          \\ \hline
\multirow{2}{*}{\textbf{5}} & 13              & \multirow{2}{*}{0.193}            & \multirow{2}{*}{0.745} & \multirow{2}{*}{403,217} \\ \cline{2-2}
                            & 39              &                                   &                        &                          \\ \hline
\multirow{2}{*}{\textbf{6}} & 12              & \multirow{2}{*}{0.222}            & \multirow{2}{*}{0.633} & \multirow{2}{*}{419,554} \\ \cline{2-2}
                            & 41              &                                   &                        &                          \\ \hline
\multirow{2}{*}{\textbf{7}} & 12              & \multirow{2}{*}{0.205}            & \multirow{2}{*}{0.692} & \multirow{2}{*}{704,289} \\ \cline{2-2}
                            & 42              &                                   &                        &                          \\ \hline
\multirow{2}{*}{\textbf{8}} & 11              & \multirow{2}{*}{0.224}            & \multirow{2}{*}{0.658} & \multirow{2}{*}{720,316} \\ \cline{2-2}
                            & 43              &                                   &                        &                          \\ \hline
\end{tabular}\label{table:cost_darcy_E}
\end{table}

Direct comparisons between our approach and existing methods in the sOED literature are not straightforward.
To evaluate the effectiveness and accuracy of our approach, we adapt two widely used batch OED methods to the sequential setting.
Specifically, we compare our method against a nested Monte Carlo estimator of the iEIG (as detailed in Supplementary Material), which is computationally expensive but serves as the gold standard for nonlinear Bayesian inverse problems.
Additionally, we compare against a method based on multiple linearizations and Gaussian approximations proposed in~\cite{WuChenGhattas:2023:1} (referred to as multi-linearization from hereinafter and detailed in Supplementary Material), which is also computationally efficient but may be inaccurate for highly non-Gaussian posteriors.

\begin{figure}
\vspace{-15pt}
\centering
\begin{tikzpicture}
\node[inner sep=0pt] (a) at (0,0)
{\includegraphics[width=0.93\textwidth, trim = {0 2em 0 3em}, clip]{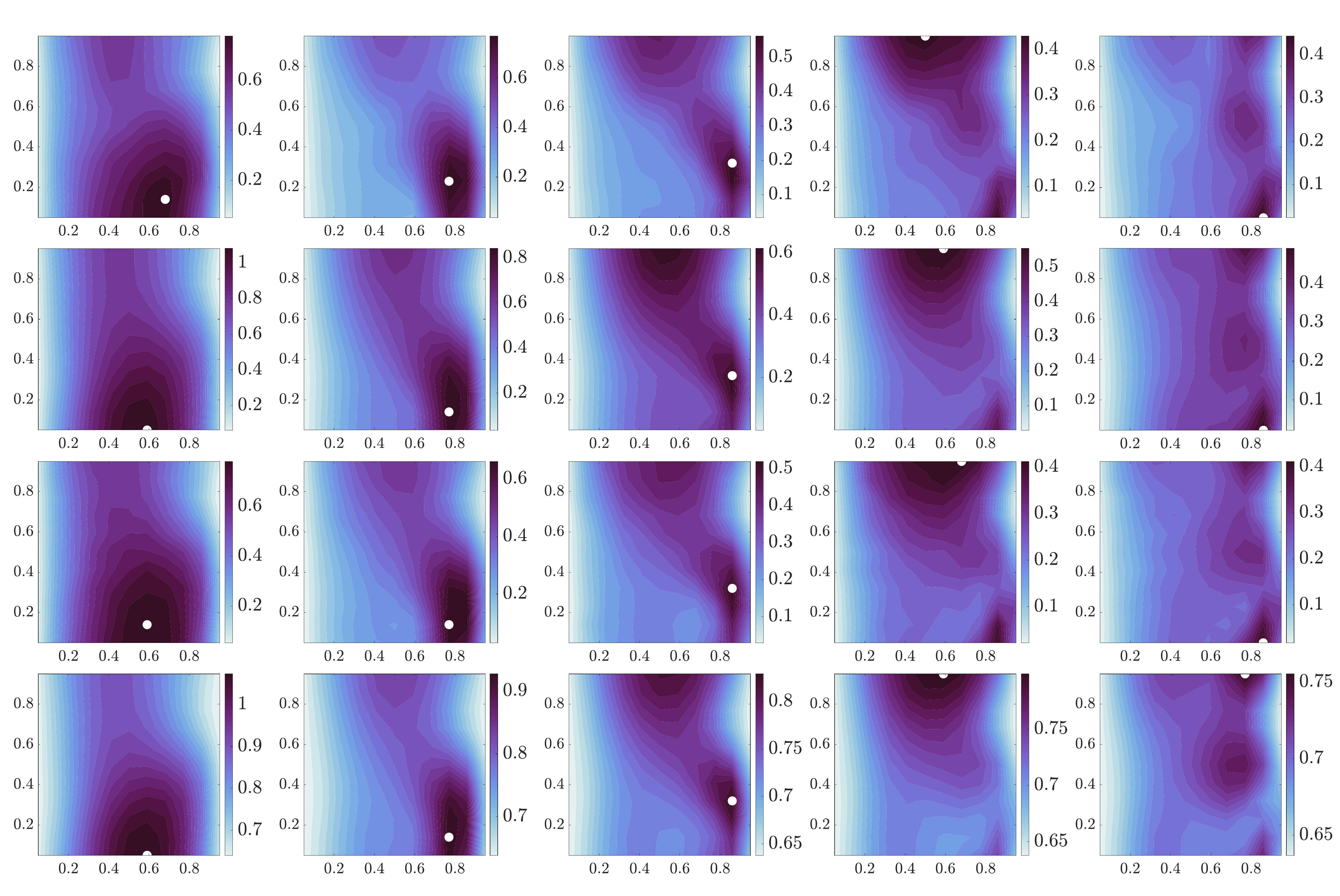}};
\node at (-7.7,3.5) {{${\EIG}^{\text{C}}$}};
\node at (-7.7,1.2) {{${\EIG}^{\text{I}}$}};
\node at (-7.7,-1.1) {{${\EIG}^{\text{NMC}}$}};
\node at (-7.7,-3.45) {{${\EIG}^{\text{GN}}$}};
\label{fig:Darcy_iEIG_comp}
\end{tikzpicture}\vspace{-20pt}
\caption{The upper bounds on the incremental EIG for stages $k = 1,\ldots,5$ (top two rows, from left to right) compared with a nested Monte Carlo estimator (third row) and a linearization-based estimator (bottom row). The upper bound was computed using~\cref{Alg:DLIS_combined} with $N = 100$ samples to approximate $\Hmat_{\mathrm{I}}^k$ at each stage. The nested Monte Carlo estimates were computed using $N = \num{10000}$ samples from the joint $\likelihood{\data_k}{\params}\,\widehat{\target}(\params \given \Hist_{k{-}1})$.}
\end{figure}

We begin by assessing the quality of the upper bounds.
Because the iEIG measures the information content of a design relative to a given prior, we fix the sequence of priors and transport maps $\{\cT^{k}\}_{k=1}^K$ to enable a meaningful comparison.
At each experimental stage $k$, with $\cT^{k{-}1}$ fixed, we evaluate both upper bounds~\eqref{eq:iEIG_uB} and~\eqref{eq:iEIG_entropy}, a reference nested Monte Carlo estimator of the iEIG (using $N=10^4$ samples) at each candidate location, as well as the iEIG approximated by the multi-linearization method.
The results of the first four experimental stages are visualized in~\cref{fig:Darcy_iEIG_comp}. They indicate that using 100 samples appears sufficient to capture the variation in the expected information gain for both of our bounds. While both bounds are sharp in this example, the covariance bound is consistently tighter. For the first three stages, both bounds and the multi-linearization give visually similar iEIG contours compared to those of nested Monte Carlo, albeit the magnitude of the estimated iEIG obtained by multi-linearization significantly differs from the others past the first stage. Starting at the fourth design, iEIG estimated by the multi-linearization begins to diverge from the others. A potential reason is that the posterior starts deviating from a Gaussian distribution as more data are observed.

\begin{figure}
\vspace{-5pt}
\centering
\begin{tikzpicture}
\node[inner sep=0pt] (a) at (-3,0)
{\includegraphics[width=0.85\textwidth, trim = {0 1em 1em 2em }, clip]{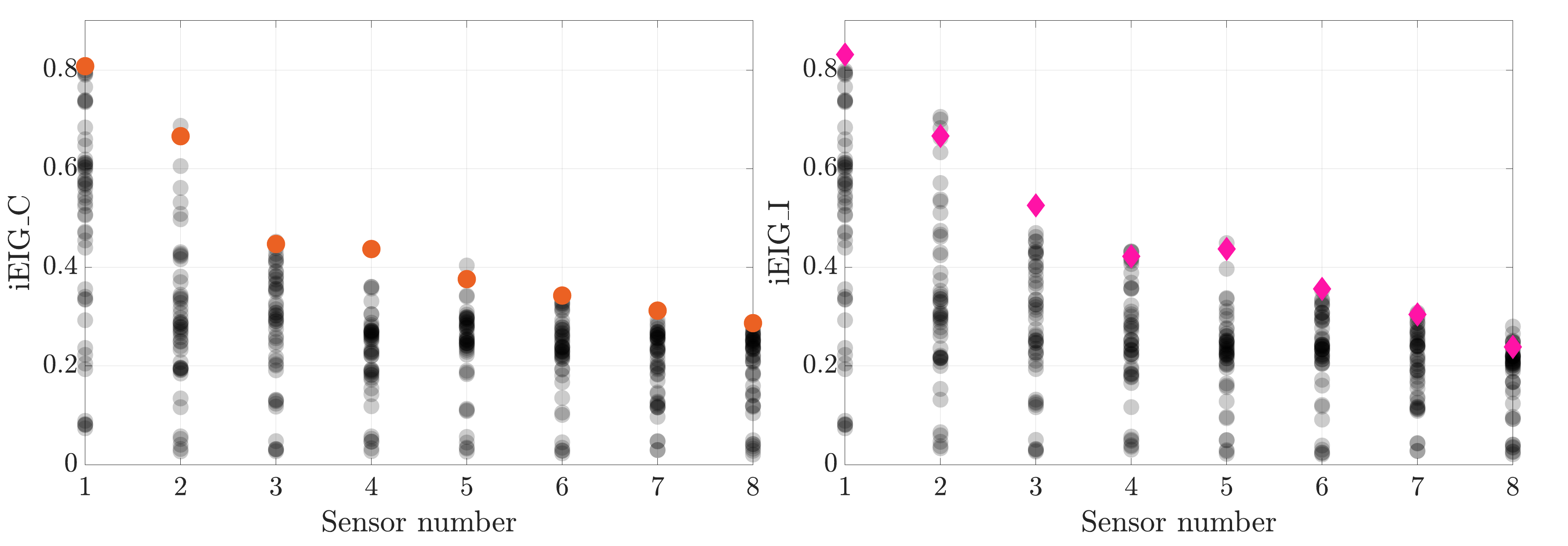}};
\end{tikzpicture}\vspace{-12pt}
\label{fig:Darcy_iEIG_compNMC}
\caption{A comparison of the iEIG for stages $1-8$ using the designs maximizing the information-based iEIG bound (right, pink diamonds) and covariance-based iEIG bound (left, orange circles) as well as 50 randomly chosen designs (black dots). The iEIG is computed using nested Monte Carlo.}
\end{figure}

In~\cref{fig:Darcy_iEIG_compNMC}, we fully evaluate the effectiveness of using the upper bounds to select sequential optimal designs. Using each bound, we run~\cref{Alg:sOED_wRestart} through the sOED cycle. At each stage $k$ and for each bound, we use nested Monte Carlo with \num{10000} samples to estimate the actual iEIG for the design obtained by the upper bound, as well as for 50 randomly selected designs. The designs obtained from those two upper bounds achieved comparable performance, and both outperform the randomly selected designs almost all the time. As more sensors are selected, the performance gap naturally decreases and is expected to shrink further with additional experiments as information accumulates.

\subsection{Problem 2: light source and pressure sensor placement for photoacoustic imaging}\label{sec:num2}

We then consider a photoacoustic imaging (PAI) application with simplified physics.
PAI is a hybrid medical imaging modality that aims to combine the high contrast of optical imaging with the high spatial resolution of ultrasound via the photoacoustic effect.
In PAI, laser-induced ultrasound waves are observed along the boundary of the tissue sample and used to reconstruct spatially varying optical properties of the tissue, such as absorption and scattering coefficients.
We refer the readers to~\cite{CoxLauferArridgeBeard:2012:1,GrohlSchellenbergDreherMaier-Hein:2021:1,LutzweilerRazansky:2013:1} for an overview of the background of PAI.
Here, we focus only on the details relevant to our example. 
We aim to infer the absorption coefficient $\mu_a$ ($\m^{-1}$), which is related to the observed data via solution of two coupled PDEs.

The first PDE represents the optical component of the PAI problem, relating tissue absorption $\mu_a$ to the initial pressure $p_0$ (Pa). Accurate modeling of light transport in tissues can be achieved using computationally intensive Monte Carlo simulations or radiative transfer equations. In highly scattering media, a common simplification is to apply the diffusion approximation to the radiative transfer equations~\cite{TarvainenCox:2024:1}.
For a two-dimensional domain $\Omega = [0,5]\times[0,3] \subset \mathbb{R}^2$, given an illumination source $\source$ on the boundary $\Gamma$ and an absorption $\mu_a$ the initial pressure $p_0$ is related to the light fluence $\phi$, which satisfies the PDE
\begin{equation}
\begin{split}
\mu_a \phi - \nabla \cdot \left( \kappa(\mu_a) \nabla \phi \right) = 0 &\quad \mbox{ for } x \in \Omega \\
\frac{1}{\target}\phi+ \frac{1}{2}\kappa(\mu_a) \nabla \phi \cdot n = \source &\quad \mbox{ for } x \text{ on } \Gamma, \\
p_0 = \Gamma \mu_a \phi &\quad \mbox{ for } x \in \Omega.
\label{eq:optic_RB}
\end{split}
\end{equation}
The diffusion coefficient $\kappa(\mu_a) = \frac{1}{2(\mu_a+\mu_s')}$ ($\m)$ depends on $\mu_a$ and the reduced scattering coefficient, $\mu_s'$ ($\m^{-1}$).
The initial pressure $p_0$ is proportional to the absorbed optical energy density $\mu_a\phi$ through the Gr\"uneisen parameter $\Gamma$, which is typically spatially-dependent and unknown in practice.
While $\Gamma$ and $\mu_s'$ could also be treated as inference parameters, doing so would significantly increase the computational cost. Moreover, simultaneous reconstruction of all three parameters ($\mu_a$, $\Gamma$, $\mu_s'$) poses theoretical challenges~\cite{BalRen:2011:1}. For our model example, we fix $\Gamma \equiv 1$ and $\mu_s' \equiv 20 \, \mathrm{cm}^{-1}$.

The second PDE is used to represent the acoustic component of the PAI problem.
Given the initial pressure $p_0$, the  pressure $p$ is obtained by solving the acoustic wave equation:
\begin{equation}
\begin{split}
\frac{1}{{c}^2}p_{tt} - \Delta p & = 0\;\, \quad \mbox{ in } \Omega \times (0,T] \\
p(t=0) &  = p_0 \quad \mbox{ in } \Omega\\
p_t(t=0) & = 0 \;\, \quad \mbox{ in } \Omega \\
\nabla p\cdot \bvec{n} - \frac{1}{c}p_t & = 0 \; \quad \mbox{ on } \Gamma\times (0,T],
\label{eq:AC}
\end{split}
\end{equation}
where the sound speed $c$ is assumed to take a constant value of $1510 \frac{m}{s}$ for the simulation.

\begin{figure}
\vspace{-15pt}
\centering
\hspace{-17pt}\begin{tikzpicture}
\node[inner sep=0pt] (a) at (0,0)
{\includegraphics[width=0.885\textwidth, trim={0 25em 0, 1em}, clip]{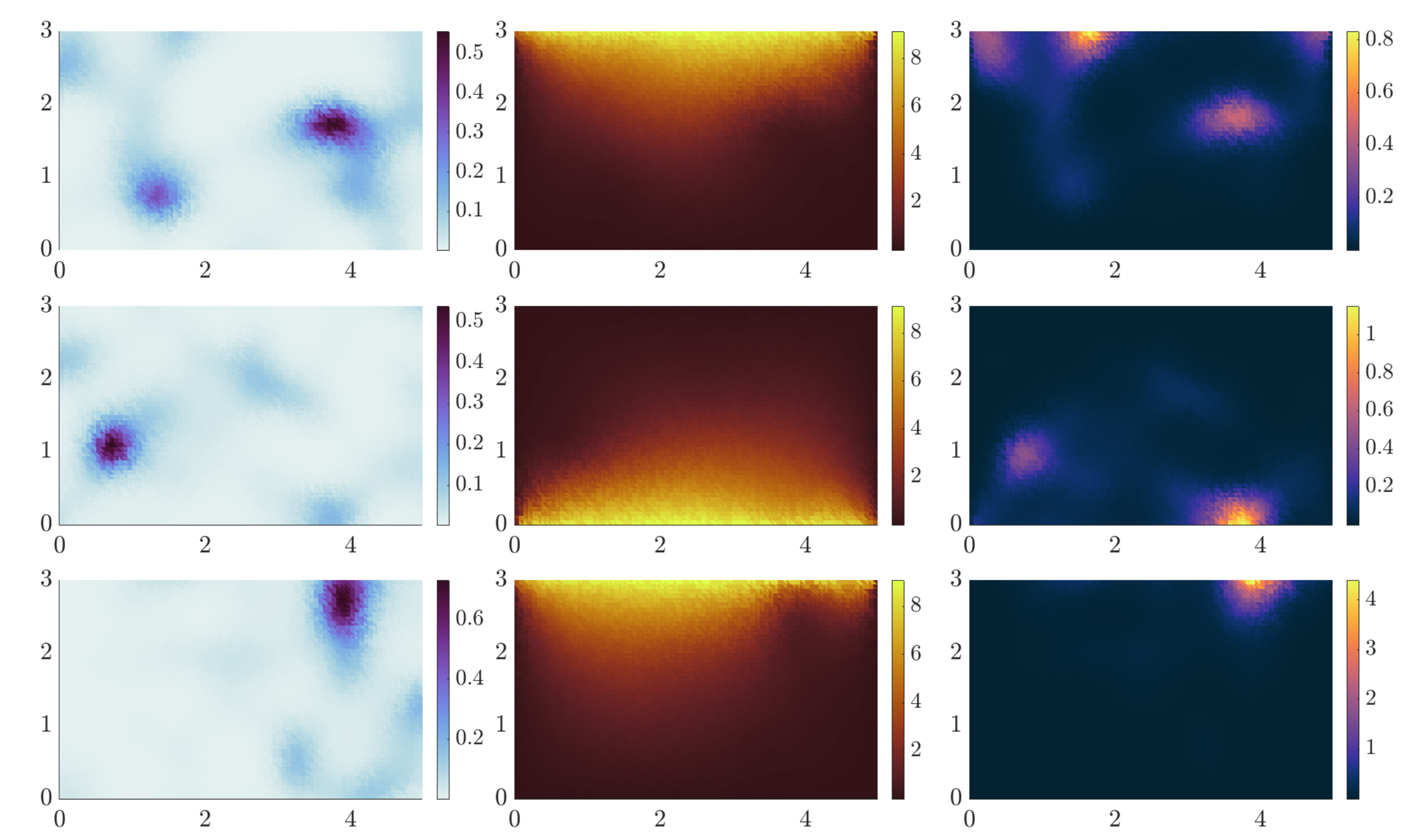}};
\node at (-4.5,2.8) {{$\mu_a$}};
\node at (0,2.8) {{$\phi(\design,\mu_a)$}};
\node at (4.5,2.8) {{$p_0(\design,\mu_a)$}};
\node at (-3.0,0.7) {{$\Omega$}};
\node at (-3.0,-2.1) {{$\Omega$}};
\node at (1.5,0.7) {{\color{white}{$\Omega$}}};
\node at (1.5,-2.1) {{\color{white}{$\Omega$}}};
\node at (6.0,0.7) {{\color{white}{$\Omega$}}};
\node at (6.0,-2.1) {{\color{white}{$\Omega$}}};
\label{fig:PAI_prSamps}
\end{tikzpicture}\vspace{-15pt}
\caption{Two sample absorption coefficients from the prior (left), the corresponding light fluence $\phi$ solving~\eqref{eq:optic_RB} (middle) with light source on the top and bottom, respectively, and the initial pressure $p_0$ (right).}
\centering
\begin{tikzpicture}
\node[inner sep=0pt] (a) at (0,0)
{\includegraphics[width=0.9\textwidth]{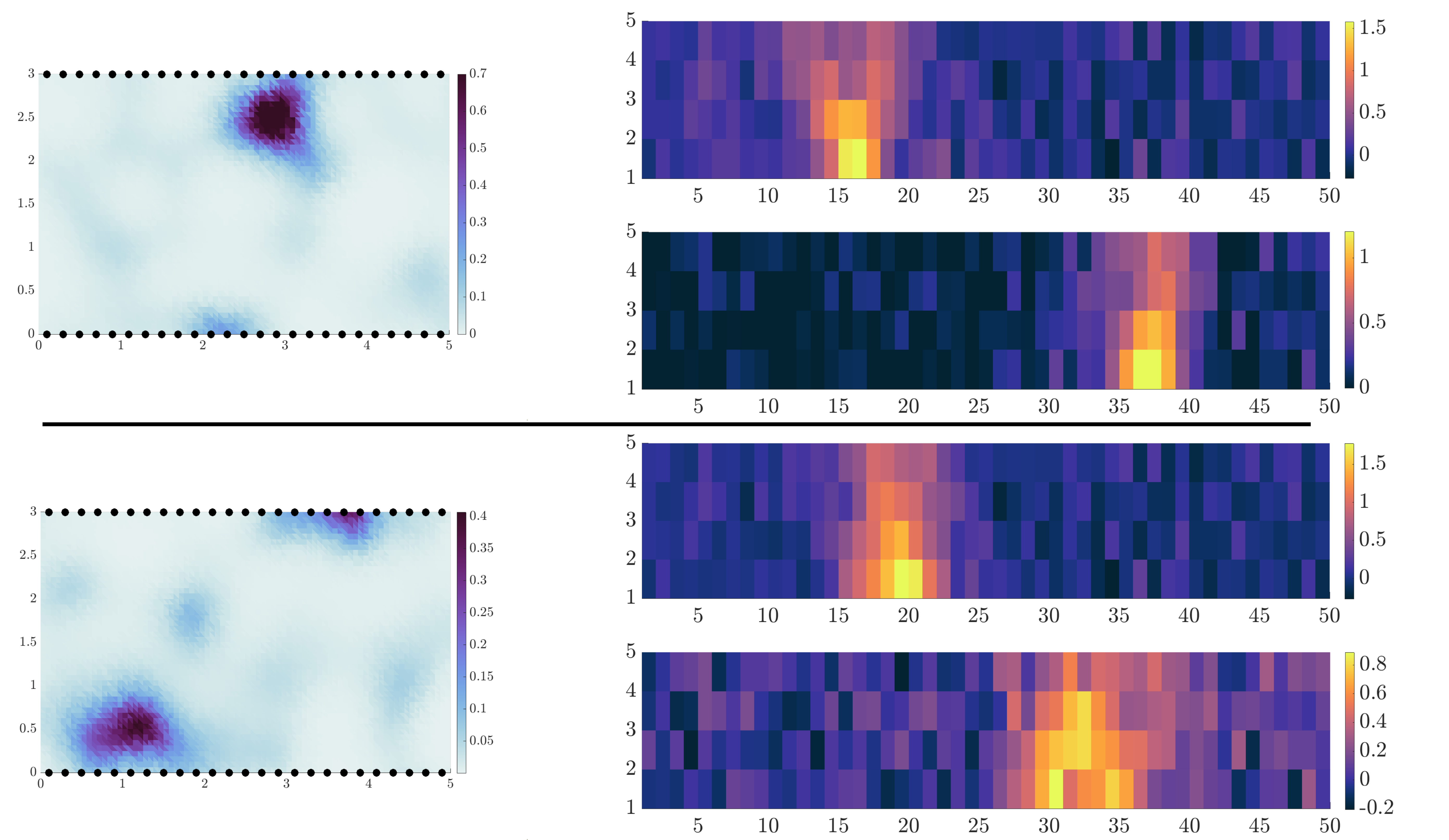}};
\node at (-4.6,3.8) {{$\mu_a^1$}};
\node at (-2.8,1.1) {{$\Omega$}};
\node at (-2.8,-3.2) {{$\Omega$}};
\node at (-4.6,-0.4) {{$\mu_a^2$}};
\node at (-1.2,3.1) {{$\tau_i$}};
\node at (-1.2,1.1) {{$\tau_i$}};
\node at (-1.2,-1.1) {{$\tau_i$}};
\node at (-1.2,-3.1) {{$\tau_i$}};
\node at (3,-4.3) {{$e_1$ (sensor number)}};
\label{fig:PAI_setup}
\end{tikzpicture}\vspace{-15pt}
\caption{Two ``true'' absorption coefficients used to synthesize data (left column). The black dots in the left column indicate the $50$ candidate locations for sensor placement. The corresponding observed pressure wave at all the candidate locations and all observation times $(\tau_i, i = 1,\ldots,5)$ is visualized in the right column, with the data in the top and bottom rows corresponding to illumination at the top and bottom boundaries,  respectively.}
\end{figure}

To ensure a nonnegative absorption coefficient and capture tissue heterogeneity, we set $\mu_a(m) = \exp(m)$ and place a Gaussian prior $\mathcal{N}(m_0,C_{\text{pr}})$ with $m_0 \equiv -4$ and $C_{\text{pr}}$ defined by the kernel $ \exp(-\frac{1}{2\ell^2} \norm{}{x-z}^{2})$ with $\ell = \frac{1}{\sqrt{5}}$.
\cref{fig:PAI_prSamps} shows two sample absorption coefficients from this prior and the corresponding light fluence and initial pressure fields.

In each experimental stage, the design involves selecting the location of the light source, which enters through the boundary condition for the optical PDE~\eqref{eq:optic_RB}, as well as the location of a single sensor where the pressure wave is measured. The sample can be illuminated from the top or bottom using a light source defined as
\begin{equation}
\source(e_1) = \left\{
     \begin{array}{@{}l@{\thinspace}l}
        3\exp\left[-\frac{1}{2}\frac{(x-2.5)^2}{25}\right] &\quad \text{if } y = e_1,\\
        0 &\quad \text{otherwise},
     \end{array}
   \right.
\end{equation}
with $e_1 \in \designSpace_1 \coloneqq \{0,3\} $.
The pressure $p$ is measured at a sensor located on the top or bottom boundary at five equally spaced observation times, starting from an initial time $\tau_0$, \ie, at $\tau_i = \tau_0+d_{\tau}(i-1)$ for $i = 1,\ldots,5$.
We assume each measurement is corrupted with mean zero noise with standard deviation $\sigma_{\eta} = 0.1$.
For simulations, we fix $\Ns = 50$ candidate locations for sensor placement (visualized in~\cref{fig:PAI_setup}) and enumerate them with an index $e_2 \in \designSpace_2 \coloneqq \{1,\ldots,50\}$.
With this setup, the sOED objective is to select a sequence of two-dimensional designs $\design = [e_1,e_2] \in \designSpace_1 \times \designSpace_2$, specifying both the laser location and the location of the pressure-reading sensor.
The effect of the light illumination location on the light fluence and initial pressure is visualized in~\cref{fig:PAI_prSamps} for different absorption coefficients.

After employing a finite element discretization with first-order Lagrange elements on a triangular mesh of size $h = \frac{1}{16}$ in each direction, the resulting discretized parameter $\params \in \mathbb{R}^{3969}$ corresponds to the nodal coefficients of the log absorption coefficient. Thus the discretized forward map, $\discPtO$ maps from the coefficient vector to the observations arising from both illumination choices at all the candidate locations, \ie, $\discPtO\colon\mathbb{R}^{3936} \rightarrow \mathbb{R}^{500}$.  The optical PDE is solved using the \texttt{FastFins} package, in which the observation operator defined by the acoustic wave equation is constructed using~\texttt{FEniCS}~\cite{AlnaesBlechtaHakeJohanssonKehletLoggRichardsonRingRognesWells:2015:1}.

\smallskip
\paragraph{Optimal designs} In this example, we compute optimal sensor placements and illumination locations for five experimental stages. To examine how the unknown parameter influences sequential designs, we selected optimal designs for two different absorption coefficients. The ``true'' absorption coefficients used to synthesize the data are visualized in~\cref{fig:PAI_setup}. For this problem, the Fisher information matrix varies significantly with the parameters, thus a larger number of samples was required to achieve a sufficiently accurate approximation of the average Fisher information matrix and the upper bound~\eqref{eq:iEIG_uB}. For our experiments, we used $N = 500$ samples with truncation tolerances $\epsilon_G = 0.01, \epsilon_I = 0.02$. As in the first example, to speed up computation of the tensor train surrogates, we build a reduced order model restricted to the data-free likelihood-informed basis at all candidate designs using the discrete empirical interpolation method. In this case, $109$ basis vectors are sufficient using a tolerance of $0.01$ and solving the ROM, which has a relative $L_2$ error of approximately 0.01, is approximately $130$ times faster than the coupled PDEs.

\begin{figure}
\vspace{-15pt}
\centering
\begin{tikzpicture}
\node[inner sep=0pt] (a) at (-3,0)
{\includegraphics[width=0.4\textwidth, trim = {0 1em 0 2em }, clip]{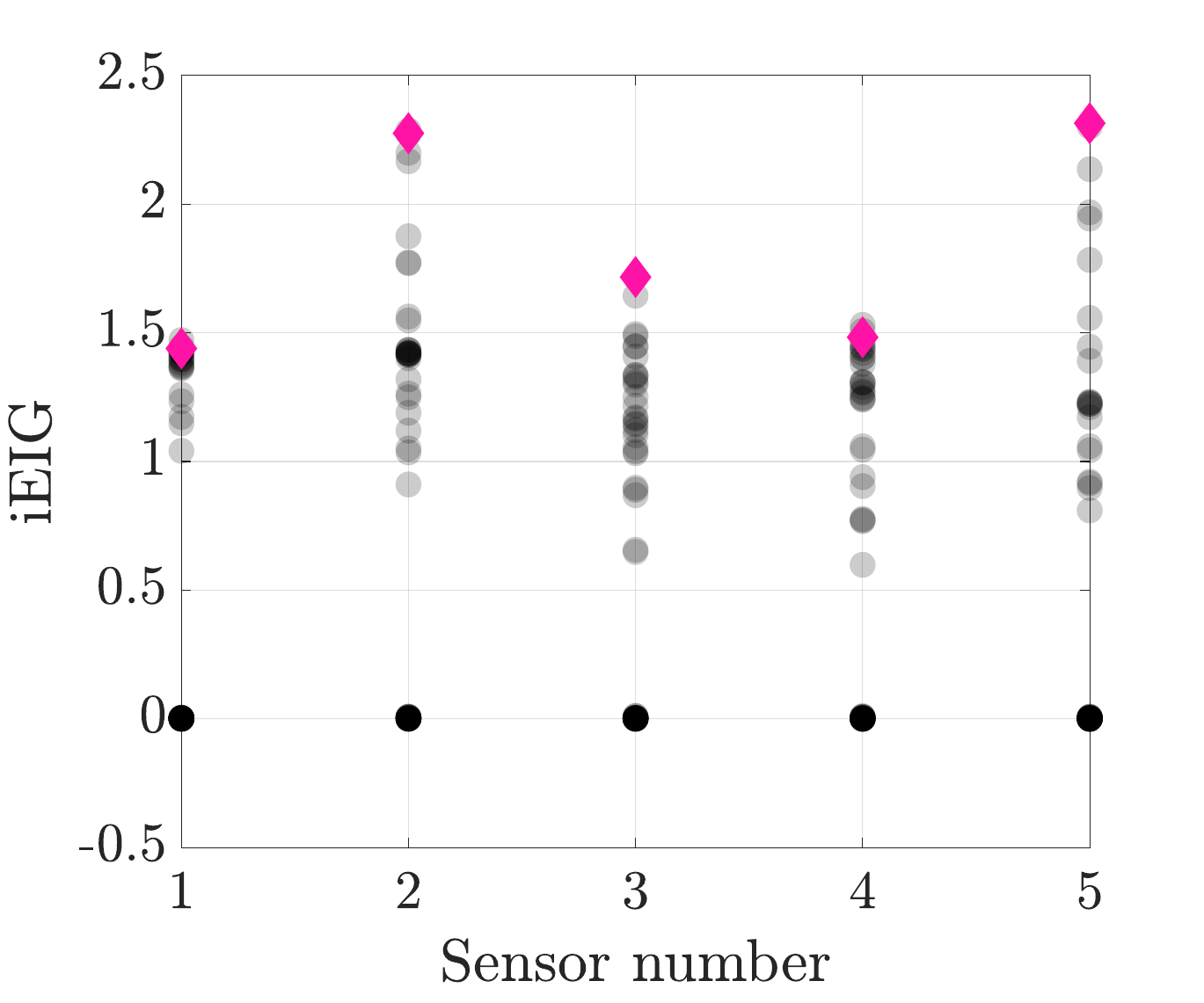}};
\node[inner sep=0pt] (b) at (4,0)
{\includegraphics[width=0.4\textwidth, trim = {0 1em 0 2em }, clip]{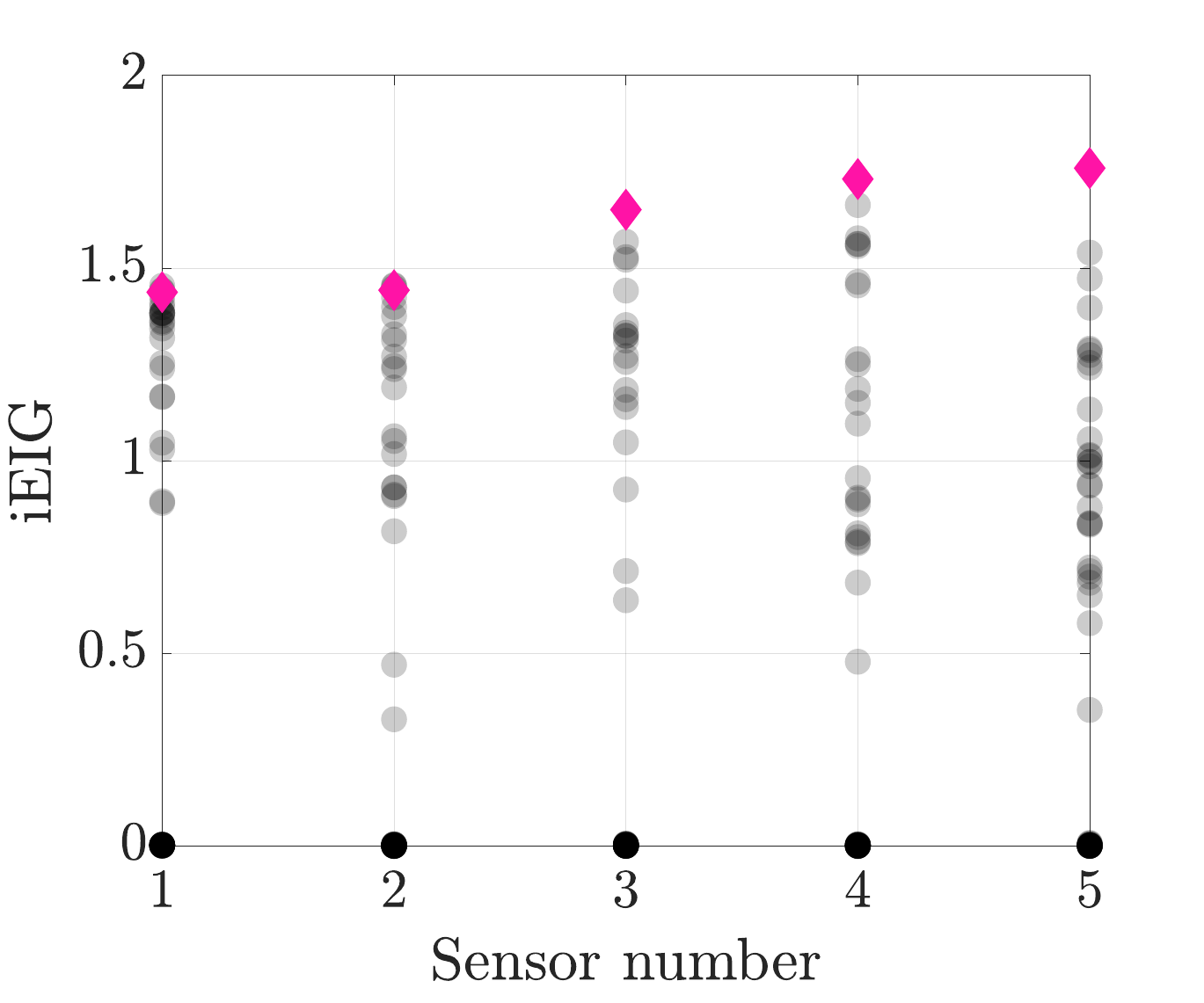}};
\label{fig:iEIG_PAI}
\end{tikzpicture}\vspace{-6pt}
\caption{A comparison of the iEIG for stages $1-5$ using the designs maximizing the iEIG upper bound (pink diamonds) as well as $50$ randomly chosen designs at each stage. The left figure corresponds to ``true'' absorption coeffient $\mu_a^1$, and the right to ``true'' absorption coefficient $\mu_a^2$.}
    \centering
    \includegraphics[width=0.9\linewidth]{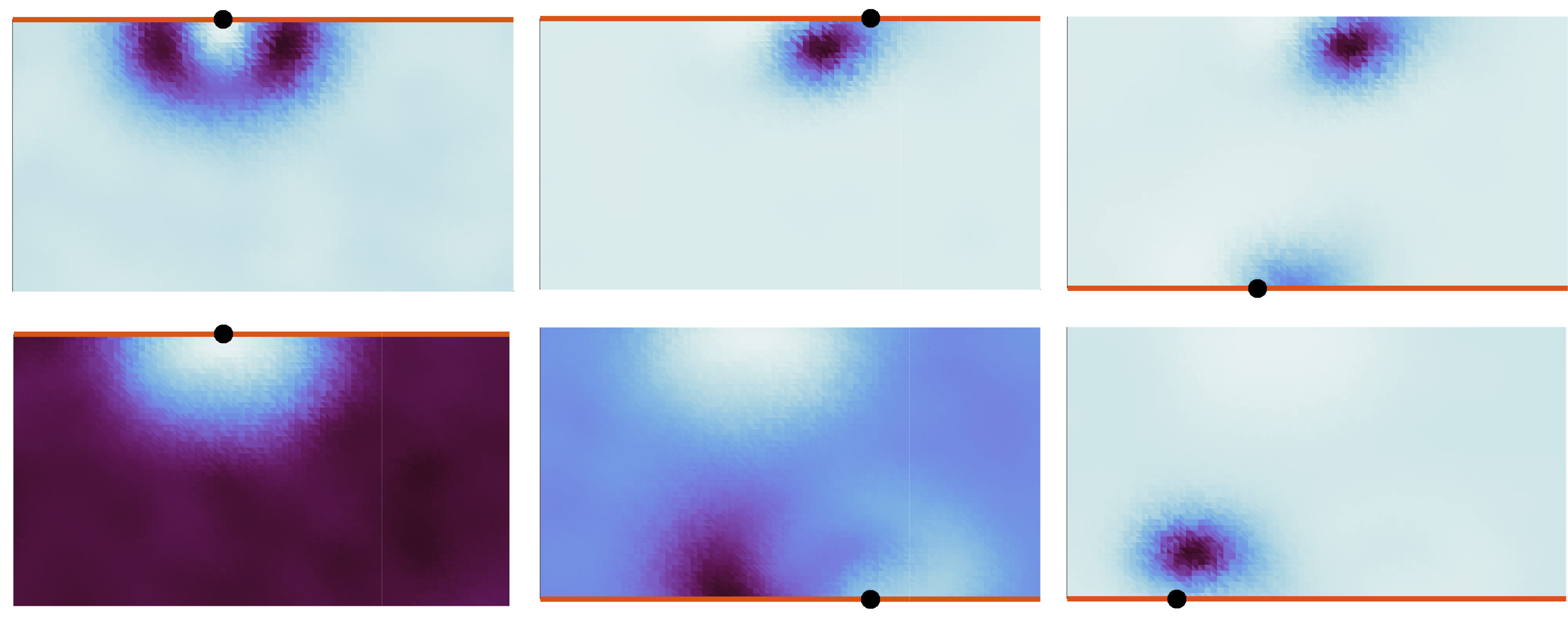}
    \caption{The posterior means for data collected at the optimal sensor locations (black dots) using the optimal illumination locations (red line) for experimental stages $1,3,5$ for both of the absorption coefficients $\mu_a^1$ (top row) and $\mu_a^2$ (bottom row) used to synthesize the data~\cref{fig:PAI_setup}. }
    \label{fig:PAI_designs}
\end{figure}

In~\cref{fig:PAI_designs}, for both absorption coefficients, we visualize the optimal laser and sensor locations for experimental stages 1-5, as well as the posterior mean corresponding to the data collected from these synthetic experiments. Our results indicate that, for our setup, it is optimal to illuminate and collect data at the same boundary. Additionally, the sequential designs appear to be highly dependent on the true parameter. For the first ``true'' absorption coefficient $\mu_a^1$, where the main inclusion is located close to the top boundary, most of the optimal designs align on the top boundary. In contrast, for the second example $\mu_a^2$, where the inclusion is near the bottom boundary, illuminations from the bottom are preferred after the initial few experiments. To further assess the effectiveness of our upper bound, we evaluate a nested Monte Carlo estimate of the iEIG at the optimal designs and compare it to the iEIG at 50 randomly chosen designs. As shown in~\cref{fig:iEIG_PAI}, our designs consistently outperform these randomly chosen designs, suggesting that the bound provides a reliable guide for optimality. In the figure, designs with an approximated iEIG of 0 correspond to configurations where the sensor is placed on the opposite boundary from the laser.
The computational cost of choosing the optimal designs and approximating the posteriors, along with the approximation errors, is presented in~\cref{table:PAI_data}.

\begin{table}[h]\label{table:PAI_data}
\footnotesize
\caption{The dimensions of data-free and data-dependent LISs (2nd column), the Hellinger error of the surrogate posterior (3rd column), ESS per sample (4th column) and the number of ROM solves for approximating the posterior (5th column) in each stage. The left table corresponds to $\mu_{a}^1$ and the right table to $\mu_a^2$.} 
\centering
\begin{tabular}{|c|c|c|c|c|}
\hline
\textbf{$k$}              & \textbf{\# LIS} & \textbf{$\mathcal{D}_{\text{H}}$} & \textbf{ESS/N}         & \textbf{\# ROMs}           \\ \hline
\multirow{2}{*}{\textbf{1}} & 9               & \multirow{2}{*}{0.152}            & \multirow{2}{*}{0.845} & \multirow{2}{*}{282,875}   \\ \cline{2-2}
                            & 17              &                                   &                        &                            \\ \hline
\multirow{2}{*}{\textbf{2}} & 14              & \multirow{2}{*}{0.161}            & \multirow{2}{*}{0.808} & \multirow{2}{*}{577,933}   \\ \cline{2-2}
                            & 20              &                                   &                        &                            \\ \hline
\multirow{2}{*}{\textbf{3}} & 13              & \multirow{2}{*}{0.169}            & \multirow{2}{*}{0.755} & \multirow{2}{*}{661,044}   \\ \cline{2-2}
                            & 26              &                                   &                        &                            \\ \hline
\multirow{2}{*}{\textbf{4}} & 9               & \multirow{2}{*}{0.394}            & \multirow{2}{*}{0.153} & \multirow{2}{*}{712,473}   \\ \cline{2-2}
                            & 46              &                                   &                        &                            \\ \hline
\multirow{2}{*}{\textbf{5}} & 19              & \multirow{2}{*}{0.318}            & \multirow{2}{*}{0.333} & \multirow{2}{*}{1,182,650} \\ \cline{2-2}
                            & 46              &                                   &                        &                            \\ \hline
\end{tabular}
\quad
\begin{tabular}{|c|c|c|c|c|}
\hline
\textbf{$k$}              & \textbf{\# LIS} & \textbf{$\mathcal{D}_{\text{H}}$} & \textbf{ESS/N}         & \textbf{\# ROMs}           \\ \hline
\multirow{2}{*}{\textbf{1}} & 9               & \multirow{2}{*}{0.023}            & \multirow{2}{*}{0.996} & \multirow{2}{*}{245,737}   \\ \cline{2-2}
                            & 10              &                                   &                        &                            \\ \hline
\multirow{2}{*}{\textbf{2}} & 9               & \multirow{2}{*}{0.135}            & \multirow{2}{*}{0.693} & \multirow{2}{*}{339,605}   \\ \cline{2-2}
                            & 25              &                                   &                        &                            \\ \hline
\multirow{2}{*}{\textbf{3}} & 16              & \multirow{2}{*}{0.139}            & \multirow{2}{*}{0.843} & \multirow{2}{*}{699,670}   \\ \cline{2-2}
                            & 33              &                                   &                        &                            \\ \hline
\multirow{2}{*}{\textbf{4}} & 17              & \multirow{2}{*}{0.285}            & \multirow{2}{*}{0.341} & \multirow{2}{*}{802,001}   \\ \cline{2-2}
                            & 39              &                                   &                        &                            \\ \hline
\multirow{2}{*}{\textbf{5}} & 17              & \multirow{2}{*}{0.233}            & \multirow{2}{*}{0.533} & \multirow{2}{*}{1,426,682} \\ \cline{2-2}
                            & 43              &                                   &                        &                            \\ \hline
\end{tabular}
\end{table}
\section{Concluding remarks}\label{sec:conclusions}
We presented a greedy sequential experimental design method for Bayesian inverse problems in high-dimensional, non-Gaussian settings. Our approach combines conditional transport maps for amortized inference with likelihood-informed subspaces for scalability. These ingredients enable efficient evaluation of sharp bounds on the incremental expected information gain using a single Monte Carlo loop, avoiding high-variance nested estimators. Numerical results for two model problems demonstrate that the bound yields informative designs in challenging practical applications.

Despite the substantial dimensionality reduction provided by the LIS, constructing transport maps remains forward-model intensive. As a result, surrogates or reduced-order models are a key component for problems with extremely costly forward models. Although the ROMs used in our experiments were sufficiently accurate, settings with non-negligible surrogate error will require appropriate model-error correction tools to mitigate bias.

Several opportunities exist to improve the proposed approach. For example, model evaluations used to construct iEIG bounds could be stored and reused (via interpolation of outputs or Jacobians) to reduce computational cost. As shown by the subsurface flow example, hybrid strategies may also be effective: Gaussian approximations to the posterior can be employed when appropriate, e.g., for preconditioning the conditional map construction, to further reduce the complexity of the construction and deployment of conditional transport maps. 

\appendix
\section{Supplementary Material} 

\subsection{Comparison of entropic and information bounds for a toy problem}\label{appendix:bound_comp}

All the bounds introduced in Section 3.1.2 and 3.1.3 are exact if both prior and posterior are Gaussian. To explore the effectiveness of the bounds outside this Gaussian setting, we consider a modified version of the example problem from~\cite[Section 5.1]{HuanMarzouk:2013:1}, 
\begin{equation}
d = e^2m^3+\frac{1}{\beta}m\exp\left(-|{0.2-e}|^{\beta}\right),\quad m \sim \mathcal{U}([0,1]), \quad \eta \sim \mathcal{N}(0,10^{-4}).
\label{eq:toy_model}
\end{equation}
Here, the observation depends nonlinearly on the design variable $e \in [0,1]$ as well as the scalar unknown parameter $m$. The term $\beta > 0$ is used to introduce a tunable peak in the EIG at $e = 0.2$, with $\beta$ controlling both the sharpness of this peak as well as the relative strength of the signal of $\frac{1}{\beta}m\exp\left(-|{0.2-e}|^{\beta}\right)$ compared to the smoother $e^2m^3$ term.

\begin{figure}
\centering
{\includegraphics[width=\textwidth]{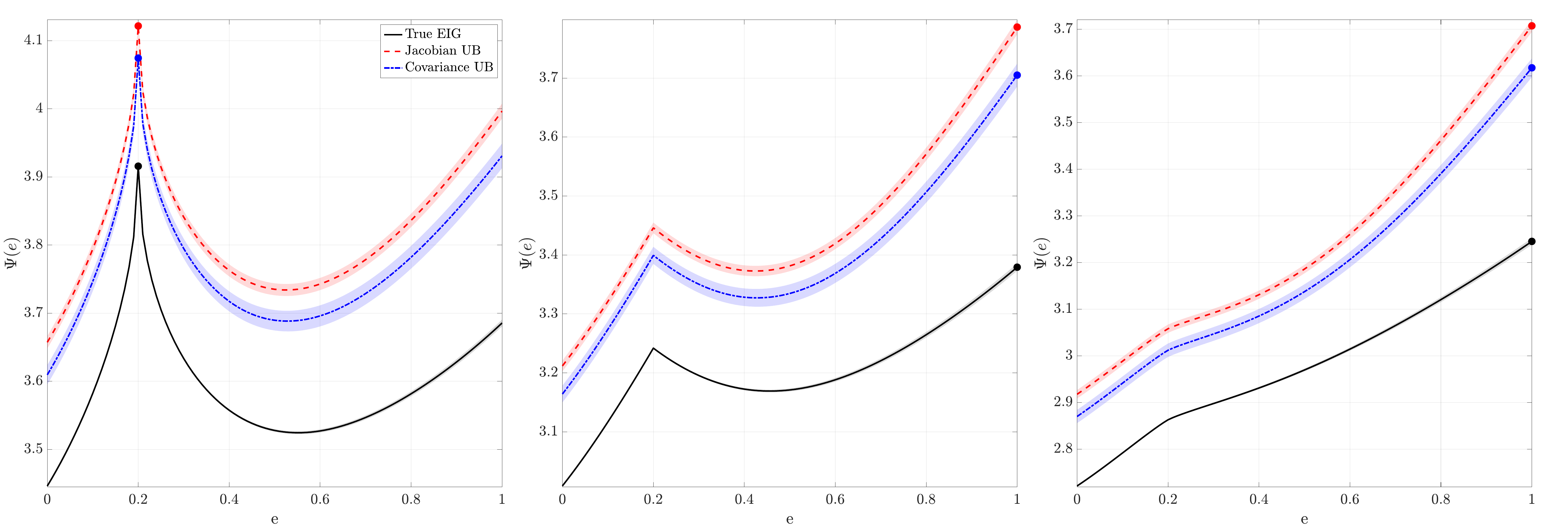}}
\label{fig:bound_comp_simple}
\caption{A comparison of the gradient-based and entropic upper bounds, $\EIG^{\mathrm{I}}$ (red), and $\EIG^{\mathrm{E}}$ (blue), respectively, along with the ``true'' EIG (black) for $e \in [0,1]$ for the quadratic-cubic model~\eqref{eq:toy_model}. The optimal design is also shown for each approximation using the filled circle.}
\end{figure}

Using the entropic form of the EIG  $\EIG(e) = h(y \given e)-\Expect{m}{h(y \given e,m)}$, the true EIG is approximated numerically using a sample-based approximation to $h(y|e)$ with $10000$ data samples combined with quadrature to determine $\pi(y \given e)$. In this case, both data and parameter are one-dimensional, so the gradient-based bounds ($\EIG^{\mathrm{G}}$ and $\EIG^{\mathrm{I}}$) are equivalent, and $\EIG^{\mathrm{C}} \leq \EIG^{\mathrm{G}} = \EIG^{\mathrm{I}}$. Both $\EIG^{\mathrm{C}}$ and $\EIG^{\mathrm{I}}$ were evaluated numerically using $N=1000$ samples. A comparison of the entropic and gradient-based bounds for $100$ runs is visualized in~\cref{fig:bound_comp_simple} for $\beta \in \{0.5,1.0,1.5\}$.
Both estimators exhibit an increasing gap from the true EIG as $e$ approaches 1, corresponding to the stronger contribution from the cubic term in~\eqref{eq:toy_model}.
As expected, the covariance-based bound is sharper than the gradient-based bound in this scalar setting, however, the gap between the two is not substantial in this example. 
The variances of the estimators are also similar, although the covariance-based estimator exhibits slightly higher variance than the information-based estimator for a fixed sample size. 
Most importantly, in this particular example, both estimators correctly identify both local maxima on average.

\subsection{Comparison method I --- Nested Monte Carlo}\label{appendix:NMC}
To evaluate the performance of our upper bound on the iEIG, we compare it with a nested Monte Carlo estimator. Note that for our likelihood model, using Bayes' law, we have $\target(\params\given\design_k,\data_k,\Hist_{k-1}) = \frac{\likelihood{\data_k}{\params,\design_k}\,\target(\params\given \Hist_{k-1})}{\target(\data_k \given \design_k,\Hist_{k-1})}$, and the incremental EIG can be written as follows: 
\begin{equation}
\begin{split}
    \EIG_k(\design_k) &= \Expect{\data_{k}\given \design_{k},\Hist_{k-1}}{\distKL{\target(\cdot \given \design_k,\data_k,\Hist_{k-1})}{\target(\cdot \given \Hist_{k-1})}} \\
    &= \int \int \log \left( \frac{\target(\params \given \design_k,\data_k,\Hist_{k-1})}{\target(\params\given \Hist_{k-1})}\right) \target(\params \given \design_k,\data_k,\Hist_{k-1}) \, \mathrm{d}\params\,  \target(\data_k \vert \design_k, \Hist_{k-1}) \, \mathrm{d}\data_k \\
    &= \int  \int \log \left( \frac{\likelihood{\data_k}{\params,\design_k}}{\target(\data_k \given \design_k, \Hist_{k-1})}\right) \likelihood{\data_k}{\params,\design_k}\, \target(\params\given \Hist_{k-1}) \, \mathrm{d}\params \, \mathrm{d}\data_k.
\end{split}
\end{equation}
Following~\cite{HuanMarzouk:2013:1}, $\EIG_k(\design_k)$ can be approximated at any design $\design_k \in \designSpace$ via the following double-loop or nested Monte Carlo estimator
\begin{equation}
    \EIG_k(\design_k) \approx \EIG_k^{\text{NMC}}(\design_k) = \frac{1}{N_{\text{out}}}\sum_{i=1}^{N_{\text{out}}} \left( \log(\likelihood{\data_k^{(i)}}{\params^{(i)},\design_k})-\log(\hat{\target}(\data_k^{(i)} \given \design_k, \Hist_{k-1})\right), 
\end{equation}
where $\params^{(i)}$ are drawn from the prior $\target(\params \given \Hist_{k-1})$ and $\data_k^{(i)}$ are drawn from the likelihood $\likelihood{\data_k}{\params^{(i)},\design_k}$. Since the evidence $\target(\data_k \given \design_k, \Hist_{k-1})$ typically does not have a closed-form expression, we estimate it using an inner Monte Carlo estimator: 
\begin{equation}
    \target(\data_k^{(i)} \given \design_k, \Hist_{k-1}) \approx \hat{\target}(\data_k^{(i)} \given \design_k, \Hist_{k-1}) = \frac{1}{N_{\text{in}}} \sum_{j=1}^{N_{\text{in}}} \likelihood{\data_k^{(i)}}{\params^{(i,j)},\design_k}, 
    \label{eq:NMC_evidence}
\end{equation}
where $\params^{(i,j)} \sim \target(\params\given \Hist_{k-1})$. 

To accelerate the nested Monte Carlo estimator, we also use transport maps, combined with importance sampling, to draw samples from $\target(\params\given \Hist_{k-1})$. Simulating data from the likelihood and evaluating its density involves the forward map, and thus computing the inner and outer loops for each candidate design $\design_k$ would require $O(N_{\text{in}}N_{\text{out}}\Nd)$ evaluations of the costly forward map, which can be prohibitively expensive. We mitigate the computational cost by reusing the posterior samples in two ways. First, we fix a sample set $\{\params^{(i)}\}_{i=1}^{N_{\text{out}}}$ and reuse it for each design $\design_k$. We additionally set $N_{\text{in}} = N_{{\text{out}}}$ and reuse these samples when approximating the evidence~\eqref{eq:NMC_evidence}. Sample reuse contributes to the bias of the nested Monte Carlo estimator, however, as stated in~\cite{HuanMarzouk:2013:1}, this effect is rather small. 

\subsection{Comparison method II --- Estimating iEIG using local Gaussian approximations} \label{appendix:Laplace}
We further evaluate our approach against an extension of the methods outlined in~\cite{WuChenGhattas:2023:1}. 
In particular, we extend the prior sample point approximation method described in section 4.4 of~\cite{WuChenGhattas:2023:1} to estimate the iEIG using successive Gauss-Newton approximations to the posterior. 
Employing a Monte Carlo approximation to the iEIG, we have
\begin{equation}
    \EIG_k(\design_k) \approx \EIG_k^{\mathrm{GN}}(\design_k) = \frac{1}{N}\sum_{i=1}^N\distKL{\target^G(\cdot \given \data_k^{(i)},\design_k,\Hist_{k-1})}{\target^G(\cdot \given \Hist_{k-1})},
\end{equation}
where $\target^G(\params \given \data_k,\design_k,\Hist_{k-1})$ and $\target^G(\params \given \Hist_{k-1})$ denote Gaussian approximations to the candidate stage-$k$ posterior and stage-$k$ prior, respectively. The data samples $\data_k^{(i)}$ for $i = 1,\ldots,N$ are obtained by drawing a sample $\params^{(i)}$ from the Gaussian approximation to the prior, and synthesizing the noisy data using the accurate forward map, $\PtO$. 

In particular, given a fixed history of observed data and experimental conditions $\Hist_{k-1}$, denote the corresponding maximum-a-posteriori (MAP) estimator as $\params_{k-1}^{\mathrm{MAP}}$. Under our assumption of i.i.d.\@ Gaussian noise $\noise \sim \mathcal{N}(\bvec{0},\sigma_{\eta}^2\I{})$, using $\Jac \PtO(\params_{k-1}^{\mathrm{MAP}})$ to denote the Jacobian of the parameter-to-observable map evaluated at the MAP estimator, the Gauss-Newton approximation to the stage-$k$ prior is, 
\begin{align*}
\target^G(\params \given \Hist_{k-1}) &\sim \mathcal{N}\left(\params_{k-1}^{\mathrm{MAP}},\mathcal{C}_{k-1}(\params_{\mathrm{MAP}})\right),  \\
\mathcal{C}_{k-1}(\params_{\mathrm{MAP}}) &= \left({\mathcal{C}^{-1}_{k-2}}+\frac{1}{\sigma_{\eta}^2}\Jac\PtO(\params_{k-1}^{\mathrm{MAP}})^\top\,\Jac\PtO(\params_{k-1}^{\mathrm{MAP}})\right)^{-1}. 
\end{align*}

In each experimental stage, a Gaussian approximation needs to be constructed for every set of sample data $\data_{k}^{(i)}$. The standard Gauss-Newton approximation centers the Gaussian around the MAP estimator. However, this requires solving at least $N$ optimization problems in each experimental stage (depending on the approach used). To avoid this, we use the approach outlined in~\cite[section 4.4]{WuChenGhattas:2023:1}. In this setting, the KL divergence between $\target^G(\cdot \given \data_k^{(i)},\design_k,\Hist_{k-1})$ and $\target^G(\cdot \given \Hist_{k-1})$ can be estimated as
\begin{multline}
    \distKL{\target^G(\cdot \given \data_k^{(i)},\design_k,\Hist_{k-1})}{\target^G(\cdot \given \Hist_{k-1})} = \\ \frac{1}{2}\sum_{j=1}^{k_i}\left(\log \left(1+\lambda_j(\Hmat^{(i)}_k(\design_k))\right) - \frac{\lambda_j(\Hmat^{(i)}_k(\design_k))}{1+\lambda_j(\Hmat^{(i)}_k(\design_k))}\right)
    + \| \params^{(i)} - \params_{k-1}^{\mathrm{MAP}} \|_{\mathcal{C}^{-1}_{k-1}}
    \label{eq:DKL_Gaussians}
\end{multline}
where $\lambda_{j}(\Hmat^{(i)}_k(\design_k))$ are the nonzero eigenvalues of the matrix $$\Hmat^{(i)}_k(\design_k) = \frac{1}{\sigma_{\eta}^2}\Jac\PtO(\design_k,\params^{(i)})\,\mathcal{C}_{k-1}\Jac\PtO(\design_k,\params^{(i)})^\top.$$ 
Thus, the resulting approximation to the stage $k$ iEIG using this approach is:
\begin{equation}
    \EIG_k^{\mathrm{GN}}(\design_k) \!=\! \frac{1}{2N}\!\sum_{i=1}^N \!\left( \sum_{j=1}^{k_i} \!\left(\log \left(1+\lambda_j(\Hmat^{(i)}_k(\design_k))\right) - \frac{\lambda_j(\Hmat^{(i)}_k(\design_k))}{1+\lambda_j(\Hmat^{(i)}_k(\design_k))}\right) + \| \params^{(i)} \!-\! \params_{k-1}^{\mathrm{MAP}} \|_{\mathcal{C}^{-1}_{k-1}} \!\right)\!. 
\end{equation}

\bibliographystyle{abbrv}
\bibliography{biblio}

\end{document}